\documentclass[aos,seceqn,nameyear,dvips]{arximspdf}
\usepackage{mathbh}
\usepackage{graphicx}

\doi{10.1214/09-AOS766}
\volume{38}
\issue{3}
\pubyear{2010}
\firstpage{1767}
\lastpage{1792}

\makeatletter

\newtheorem{theorem}{Theorem}
\newtheorem{cor}[theorem]{Corollary}

\newproclaim{Step}{Step}
\newproclaim{Remark}{Remark}

\newcommand{\dhhat}{{\widehat{\widehat h}}}
\newcommand{\mathbbm}{\mathbh}

\makeatother

\begin{document}
\begin{frontmatter}

\title{Asymptotics and optimal bandwidth selection for
highest density region estimation\protect\thanksref{T1}}
\runtitle{Bandwidth selection for highest density region estimation}

\thankstext{T1}{Supported in part by Australian Research Council
Grant DP055651.}

\begin{aug}
\author[A]{\fnms{R. J.} \snm{Samworth}\ead[label=e1]{r.samworth@statslab.cam.ac.uk}} and
\author[B]{\fnms{M. P.} \snm{Wand}\corref{}\ead[label=e2]{mwand@uow.edu.au}}
\runauthor{R. J. Samworth and M. P. Wand}
\affiliation{University of Cambridge and University of Wollongong}
\address[A]{Statistics Laboratory and Department\\
\quad of Pure Mathematics\\
\quad and Mathematical Statistics\\
University of Cambridge\\
Cambridge CB2 1TP\\
United Kingdom\\
\printead{e1}} 
\address[B]{Centre for Statistical\\
\quad and Survey Methodology\\
School of Mathematics\\
\quad and Applied Statistics\\
University of Wollongong\\
Northfields Avenue\\
Wollongong 2522\\
Australia\\
\printead{e2}}
\end{aug}

\received{\smonth{8} \syear{2009}}
\revised{\smonth{11} \syear{2009}}

%
\begin{abstract}
We study kernel estimation of highest-density regions (HDR). Our
main contributions are two-fold. First, we derive a uniform-in-bandwidth
asymptotic approximation to a risk that is appropriate for HDR estimation.
This approximation is then used to derive a bandwidth selection rule
for HDR estimation possessing attractive asymptotic properties.
We also present the results of numerical studies that illustrate
the benefits of our theory and methodology.
\end{abstract}

%

\end{frontmatter}

\section{Introduction}\label{sec:intro}

A highest-density region (HDR) for a measurement of interest is a
region where the underlying density function exceeds some nominal
threshold. Given a random sample from that density, HDR estimation
typically involves determination of regions where an estimated density
is high. Kernel density estimation is the most common approach, but its
performance is heavily dependent on the choice of the bandwidth
parameter. Automatic selection of a good bandwidth for HDR estimation
is the overarching goal of this article.

\begin{figure}

\includegraphics{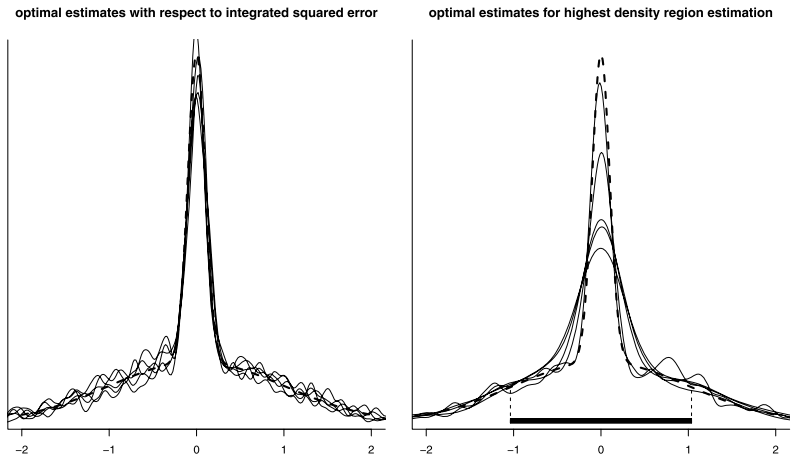}

\caption{Left panel: five kernel density estimates based
on random samples of size 1000 simulated from the density depicted by the
dashed curve. Each estimate is based on the optimal bandwidth
with respect to integrated squared error. Right panel: same as the left panel
except that the bandwidth is chosen to minimize the error for estimation
of the 20\% highest-density region. This region is shown as a thick
horizontal line at the base of the plot and its boundaries are shown
as dashed vertical lines.}
\label{fig:HDRvISE}
\end{figure}

Figure \ref{fig:HDRvISE} illustrates the bandwidth selection issue for
HDR estimation. The left panel shows five kernel density estimates
based on random samples of size 1000 from the normal mixture
$\frac{2}{3}N(0,1)+\frac{1}{3}N(0,\frac{1}{100})$ density [Density 4 of
Marron and Wand (\citeyear{MW92})]. In each case the bandwidth is
chosen to minimize the integrated squared error (ISE). In the right
panel the same random samples are used, but, instead, the bandwidths
are chosen to minimize an error appropriate for estimation of the 20\%
HDR (defined formally in Section \ref{sec:asyRisk}). This region is
shown as a thick horizontal line at the base of the plot. It is clear
from Figure \ref{fig:HDRvISE} that optimality for HDR estimation is
quite different from ISE-optimality. Low ISE requires that the two
curves be close to each other over the whole real line. However, good
estimation of the 20\% HDR only requires that the 20\% HDRs of the
kernel density estimates are close to the true region. In particular,
the sharp mode of the underlying density has no bearing upon the HDR
and there is no need to estimate it well. For this density it is
apparent that a bandwidth considerably larger than ISE-optimal
bandwidth is best for estimation of the 20\% HDR.

%

In this article we study an asymptotic risk associated with
kernel-based HDR estimation and use our theory to develop a plug-in
type bandwidth selector. Attractive asymptotic properties of our
bandwidth selector are established and good performance is illustrated
on simulated data. A self-contained function for use in the \texttt{R}
environment [R Development Core Team (\citeyear{RDevelopmentCoreTeam08})] is made available on the
Internet.

The HDR estimation problem has an established literature.
Contributions include Hartigan (\citeyear{H87}), M\"uller and Sawitzki (\citeyear{MS91}),
Polonik (\citeyear{P95}), Hyndman (\citeyear{H96}), Tsybakov (\citeyear{T97}), Ba\'illo,
Cuesta-Albertos and Cuevas (\citeyear{BCaC01}), Ba\'illo (\citeyear{B03}), Cadre (\citeyear{C06}), Jang
(\citeyear{J06}), Rigollet and Vert (\citeyear{RV09}) and Mason and Polonik (\citeyear{MP09}). Mason
and Polonik (\citeyear{MP09}) provide a thorough literature review for the
problem. Alternative terminology includes estimation of the
\textit{density contours}, \textit{density level sets} and \textit{excess mass
regions}. This literature is, however, mainly concerned with
theoretical results unconnected with the bandwidth selection problem.
Jang (\citeyear{J06}) is an applied paper on the use of HDR estimation for
astronomical sky surveys. However, the bandwidths used there are
chosen via classical ISE-based plug-in strategies. The present paper
is, to our knowledge, the first to derive theory and bandwidth
selection rules that are specifically tailored to the HDR estimation
problem.

While our proposed practical bandwidth selector relies on asymptotic
approximations, its development comes at a time when sample sizes in
applications that benefit from smoothing techniques are becoming very
large. The area of application that led to this research, flow
cytometry, typically has sample sizes in the hundreds of thousands.
The astronomical application in Jang (\citeyear{J06}) involves sample sizes in
the tens of thousands. Another HDR application is approximation of
the highest posterior density region of a parameter in a Bayesian
analysis, where only a sample from that density is available. In this
situation, the sample, most typically obtained using Markov chain Monte
Carlo methods, can arbitrarily large in size.

Section \ref{sec:asyRisk} presents an approximation to the HDR
asymptotic risk. Numerical studies support its use for
bandwidth selection. In Section \ref{sec:bwSel} we
describe plug-in strategies for bandwidth selection.
Asymptotic performance results are established and a simulation
study demonstrates practical efficacy. We conclude
with an example on daily temperature maxima in
Melbourne, Australia. Proofs are deferred
to an \hyperref[app]{Appendix}.

\section{Asymptotic risk results}\label{sec:asyRisk}

Let $f$ be a probability density function
on the real line. For $\tau \in(0,1)$, define
\[
f_{\tau } = f_{\tau }(f) = \inf\biggl\{y \in(0,\infty)\dvtx
\int_{-\infty}^\infty f(x) \mathbbm{1}_{\{f(x) \geq y\}}\, dx \leq
1-\tau \biggr\}.
\]
We call $R_\tau = \{x \in\mathbb{R}\dvtx f(x) \geq f_{\tau }\}$ the
$100(1-\tau )$\% \textit{highest-density region} of $f$
[cf. Hyndman (\citeyear{H96})]. If $(X_n)$ is
a sequence of independent random variables with density $f$, the
kernel estimator of $f(x)$ based on $X_1,\ldots,X_n$ is
\[
\widehat{f}_h(x) = \frac{1}{nh}\sum_{i=1}^n K \biggl(\frac{x-X_i}{h} \biggr),
\]
where $K\dvtx\mathbb{R} \rightarrow\mathbb{R}$ satisfies $\int K(x)\, dx
= 1$, and is called a \textit{kernel} and $h > 0$ is called the
\textit{bandwidth}. Let $\widehat{f}_{h,\tau } = f_\tau (\widehat{f}_h)$
denote the plug-in estimator of $f_{\tau }$, so that
\[
\widehat{f}_{h,\tau } = \inf\biggl\{y \in(0,\infty)\dvtx
\int_{-\infty}^\infty
\widehat{f}_h(x) \mathbbm{1}_{\{\widehat{f}_h(x) \geq y\}} \,dx \leq
1-\tau
\biggr\}.
\]
The corresponding plug-in estimator of $R_\tau $ is then
$\widehat{R}_{h,\tau } = \{x \in\mathbb{R}\dvtx \widehat{f}_h(x) \geq
\widehat{f}_{h,\tau }\}$.

Given two Borel subsets $A$ and $B$ of $\mathbb{R}$, we define their
proximity through a measure on their symmetric difference $A \triangle
B = (A \cap B^c) \cup(A^c \cap B)$. The particular measure $\mu_f$
we consider is given by
\[
\mu_f(C) = \int_C f(x) \,dx
\]
for all Borel subsets $C$ of $\mathbb{R}$. The error
$\mu_f(\widehat{R}_{h,\tau } \triangle R_\tau )$ is then then the
probability of an observation from $f$ lying in precisely one of
$\widehat{R}_{h,\tau }$ and $R_\tau $. Compared with Lebesgue
measure, $\mu_f$ puts more weight on regions where the data will tend
to be denser. It also has the advantage of admitting a simple Monte
Carlo approximation. This is important in higher-dimensional settings
where exact computation of $\mu_f(C)$ is difficult.

In Theorem \ref{Thm:Mainthm}, we derive a uniform-in-bandwidth
asymptotic expansion for the \textit{risk}
$\mathbb{E}\{\mu_f(\widehat{R}_{h,\tau } \triangle R_\tau )\}
$, which
can facilitate a theoretical, optimal choice of bandwidth
(cf. Corollary \ref{Cor:Opth}). This in turn motivates practical
bandwidth selection algorithms whose performance is studied in
Theorems \ref{Thm:Hath} and \ref{Thm:HathHO}. We will make use of the
following conditions on the underlying density, bandwidth sequence and
kernel:
\begin{enumerate}[(A1):]
\item[(A1):] $f$ is uniformly continuous on $\mathbb{R}$.
There exist finitely many points $x_1 < \cdots< x_{2r}$ such that
$f(x_j) = f_\tau $ for $j=1,\ldots,2r$, and moreover there exists
$\delta> 0$ such that $f$ is twice continuously differentiable in
$\bigcup_{j=1}^r [x_{2j-1}-\delta,x_{2j}+\delta]$ with $f'(x_{2j-1}) > 0$
and $f'(x_{2j}) < 0$ for $j = 1,\ldots,r$.
\item[(A2):] Let $h^- = h_n^-$ and $h^+=h_n^+$ be
nonnegative sequences such that $h^- \leq h^+$, such that
$n(h^-)^4/\sqrt{\log(1/h^-)} \rightarrow\infty$ and such that $h^+
\rightarrow0$ as $n \rightarrow\infty$. Then $h = h_n$ is a
sequence with $h_n^- \leq h_n \leq h_n^+$ for all $n$.
\item[(A3):] The kernel $K$ is nonnegative, continuously
differentiable, of bounded variation, and satisfies $\int x K(x)\,
dx = 0$ and $\mu_2(K) \equiv\int x^2 K(x) \,dx < \infty$.
Moreover, $K'$ is of bounded variation, and satisfies $\int
K'(x)^2 \,dx < \infty$.
\end{enumerate}
Assumption (A1) in particular implies that $f$ has a $\gamma
$-exponent
with $\gamma=1$ at level $f_{\tau }$---in other words, there
exists $C > 0$ such that
\[
\mu_f(\{x \in\mathbb{R}\dvtx |f(x) - f_{\tau }| \leq\varepsilon\})
\leq
C\varepsilon
\]
for sufficiently small $\varepsilon> 0$. This type of assumption is
common in the literature for this problem [cf. Polonik (\citeyear{P95}),
Rigollet and Vert (\citeyear{RV09})]. Although there are many parts to
condition (A3), none is very restrictive. Under
(A1), $f_{\tau }$ is the unique positive real number
satisfying $\int f(x) \mathbbm{1}_{\{f(x) \geq f_{\tau }\}} \,dx =
1-\tau $. In fact, in the course of the proof of
Theorem \ref{Thm:Mainthm} below, we will show that under
conditions (A1), (A2) and (A3),
$\widehat{f}_{h,\tau }$ has an analogous property: that is, with
probability one, for all $n$ sufficiently large,
$\widehat{f}_{h,\tau }$ is the unique positive real number satisfying
\[
\int\widehat{f}_h(x)
\mathbbm{1}_{\{\widehat{f}_h(x) \geq\widehat{f}_{h,\tau }\}} \,dx =
1-\tau .
\]

Let $\Phi$ and $\phi$ denote the standard normal distribution function
and density function, respectively, and write $R(K) = \int K^2(x)\,
dx$. Define the quantities
%
%
\begin{eqnarray}\qquad\quad
\label{Eq:D1D2D3}
D_1 &=& \frac{1}{2}\mu_2(K) \Biggl\{\sum_{j=1}^{2r} \frac{1}{|f'(x_j)|} \Biggr\}^{-1}
\Biggl[\sum_{j=1}^{2r} \frac{f''(x_j)}{|f'(x_j)|}\nonumber\\[-1pt]
&&\hspace*{122pt}{} +
\frac{1}{f_{\tau }}\sum_{j=1}^r \{f'(x_{2j}) - f'(x_{2j-1})\} \Biggr],
\nonumber\\[-9pt]\\[-9pt]
D_2 &=& R(K)f_\tau \Biggl\{\sum_{j=1}^{2r} \frac{1}{|f'(x_j)|} \Biggr\}^{-2}
\sum_{j=1}^{2r} \frac{1}{f'(x_j)^2} \quad\mbox{and} \nonumber\\[-1pt]
D_{3,j} &=& \frac{R(K)f_{\tau }}{|f'(x_j)|} \Biggl\{\sum_{k=1}^{2r}
\frac{1}{|f'(x_k)|} \Biggr\}^{-1},\qquad j=1,\ldots,2r.\nonumber
\end{eqnarray}
\begin{theorem}
\label{Thm:Mainthm}
Assume \textup{(A1)}, \textup{(A2)} and \textup{(A3)}. Then
\begin{eqnarray*}
\mathbb{E}\{\mu_f(\widehat{R}_{h,\tau } \triangle R_\tau )\} &=&
\sum_{j=1}^{2r} \biggl[\frac{B_{1,j}\phi(B_{2,j} n^{1/2}h^{5/2})}{(nh)^{1/2}}
\\
&&\hspace*{19.2pt}{} + B_{3,j} h^2
\{2\Phi(B_{2,j} n^{1/2}h^{5/2}) - 1\} \biggr]+ o \biggl(\frac{1}{(nh)^{1/2}} +
h^2 \biggr)
\end{eqnarray*}
as $n \rightarrow\infty$, uniformly for $h \in[h^-,h^+]$, where
\begin{eqnarray*}
B_{1,j} &=& 2f_{\tau } \frac{\{R(K)f_{\tau } -2D_{3,j}
+ D_2\}^{1/2}}{|f'(x_j)|},\\
B_{2,j} &=& \frac{|{1/2}\mu_2(K)f''(x_{j}) - D_1|}
{\{R(K)f_{\tau } -2D_{3,j} + D_2\}^{1/2}}
\end{eqnarray*}
and
\[
B_{3,j} = f_{\tau } \frac{|{1/2}\mu_2(K)
f''(x_{j}) - D_1|}{|f'(x_j)|}.
\]
\end{theorem}

The nature of this result is somewhat different from the results in
the existing literature which have tended to focus (sometimes in more
general settings) on the order in probability or almost surely of
$\mu_f(\widehat{R}_{h,\tau } \triangle R_\tau )$ or related
measures [e.g., Ba\'illo, Cuesta-Albertos and Cuevas (\citeyear{BCaC01}), Ba\'illo
(\citeyear{B03})]. More recent works have derived results on the limiting
behavior of suitably scaled and/or centered versions of
$\mu_f(\widehat{R}_{h,\tau } \triangle R_\tau )$ [e.g., Cadre
(\citeyear{C06}), Mason and Polonik (\citeyear{MP09})]. Rigollet and Vert (\citeyear{RV09}) provide a
finite sample upper bound for the risk, uniformly over certain
H\"older classes, with an unspecified constant in the bound. While
these theoretical results are certainly of considerable interest, our
aim in providing the asymptotic expansion in Theorem \ref{Thm:Mainthm}
is to facilitate practical bandwidth selection algorithms for this
problem---see Section \ref{sec:bwSel}.

In the course of the proof of Theorem \ref{Thm:Mainthm}, it is shown
that
\[
R(K)f_{\tau } - 2D_{3,j} + D_2 = \lim_{n \rightarrow\infty}
(nh)\operatorname{Var}\bigl(\widehat{f}_h(x_j) - \widehat{f}_{h,\tau }\bigr) > 0,
\]
so that each $B_{1,j}$ is positive. Moreover $B_{2,j}$ and $B_{3,j}$
are nonnegative, and are positive for at least one $j$. Indeed,
$B_{2,j}$ and $B_{3,j}$ are certainly positive whenever $f''(x_j) \geq
\sum_{k=1}^{2r} w_k f''(x_k)$, where the weights $w_k \propto
1/|f'(x_k)|$ sum to 1. However, this condition on $f''(x_j)$ is far
from necessary for $B_{2,j}$ and $B_{3,j}$ to be positive.

It is easily seen from Theorem \ref{Thm:Mainthm} that for any sequence
of bandwidths satisfying~(A2), if $nh^5$ is not bounded away
from zero and infinity then
$n^{2/5}\mathbb{E}\{\mu_f(\widehat{R}_{h,\tau }\times \triangle
R_\tau )\} \rightarrow\infty$ along a subsequence. On the other
hand, if $nh^5$ is bounded away from zero and infinity, then
$n^{2/5}\mathbb{E}\{\mu_f(\widehat{R}_{h,\tau }\triangle
R_\tau )\}$ is bounded. Notice that all such sequences are
permitted by the condition (A2). Focusing our attention on
bandwidth sequences of order $n^{-1/5}$ and substituting $x =
n^{1/2}h^{5/2}$, we have
\[
\lim_{n \rightarrow\infty} n^{2/5}\mathbb{E}\{\mu_f(\widehat
{R}_{h,\tau } \triangle
R_\tau )\} = \sum_{j=1}^{2r}
\biggl[\frac{B_{1,j}\phi(B_{2,j} x)}{x^{1/5}} +
B_{3,j} x^{4/5} \{2\Phi(B_{2,j} x) - 1\} \biggr].
\]
Writing this limit as $g(x) \equiv\sum_{j=1}^{2r} g_j(x)$, we see
that $g$ is continuous on $(0,\infty)$ with $g(x) \rightarrow\infty$
as $x \rightarrow0^+$ and as $x \rightarrow\infty$, so $g$ attains
its minimum. If $j$ is such that $B_{2,j}$ and $B_{3,j}$ are
positive, then it can be shown (cf. the proof of
Corollary~\ref{Cor:Opth} below), that $g_j$ has a unique minimum.
This unique minimizer represents the asymptotically optimal bandwidth
for estimating the risk in a small neighborhood of~$x_j$. Although
we typically expect the minimum of $g$ to be unique, the complicated
nature of the function $g$ and the coefficients $B_{1,j}$, $B_{2,j}$
and $B_{3,j}$ make it difficult to prove this assertion without
additional conditions. The following corollary gives the desired
result in one restricted case; however, we anticipate that the result
in fact holds much more widely.
\begin{cor}
\label{Cor:Opth}
Assume \textup{(A1)} and \textup{(A3)}. Assume further that in
\textup{(A1)} we have $r=1$ and the underlying density $f$ is
symmetric about some point on the real line. Then there exists a
unique $c_{\mathrm{opt}} \in(0,\infty)$, depending on $f$ and $K$ but
not $n$, such that any sequence of bandwidths $(h_{\mathrm{opt}})$
that minimizes $\mathbb{E}\{\mu_f(\widehat{R}_{h,\tau } \triangle
R_\tau )\}$ satisfies
\[
h_{\mathrm{opt}} = c_{\mathrm{opt}}n^{-1/5}\{1 + o(1)\}
\]
as $n \rightarrow\infty$.
\end{cor}

The additional hypotheses on $f$ imply that $B_{1,j}$, $B_{2,j}$ and
$B_{3,j}$ do not depend on $j$, and in fact in the presence of
(A1) and (A3), the conclusion of the corollary also
holds under this (weaker) condition, as can be seen from the proof.

\subsection{Numerical assessment of risk approximation}

Theorem \ref{Thm:Mainthm} yields the asymptotic risk approximation
%
%
\begin{eqnarray}\label{eq:asyRiskSim}
\mathbb{E}\{\mu_f(\widehat{R}_{h,\tau } \triangle R_\tau )\}
&\simeq&
\sum_{j=1}^{2r} \biggl[\frac{B_{1,j}\phi(B_{2,j} n^{1/2}h^{5/2})}{(nh)^{1/2}}
\nonumber\\[-8pt]\\[-8pt]
&&\hspace*{18.8pt}{} + B_{3,j} h^2
\{2\Phi(B_{2,j} n^{1/2}h^{5/2}) - 1\} \biggr].\nonumber
\end{eqnarray}

In Section \ref{sec:bwSel} we use the right-hand side of
(\ref{eq:asyRiskSim}) to develop plug-in bandwidth selection
strategies. However, it is prudent to first assess the
quality of this approximation to the risk.
We now do this through some numerical examples.

For a given $f$, $h$ and $\tau $, the risk
$\mathbb{E}\{\mu_f(\widehat{R}_{h,\tau } \triangle R_\tau )\}$
is very difficult to obtain exactly. Instead, we work with
a Monte Carlo approximation,
%
%
\begin{equation}\label{eq:riskMC}
\frac{1}{M}\sum_{i=1}^M \mu_f\bigl(\widehat{R}_{h,\tau }^{[i]}
\triangle
R_\tau \bigr),
\end{equation}
where $\widehat{R}_{h,\tau }^{[1]},\ldots,\widehat{R}_{h,\tau }^{[M]}$
are $M$ simulated realisations of $\widehat{R}_{h,\tau }$.
For large $M$ (\ref{eq:riskMC}) serves as reasonable
proxy for $\mathbb{E}\{\mu_f(\widehat{R}_{h,\tau } \triangle
R_\tau
)\}$
and is henceforth referred to as the ``exact'' risk.

%
\begin{figure}

\includegraphics{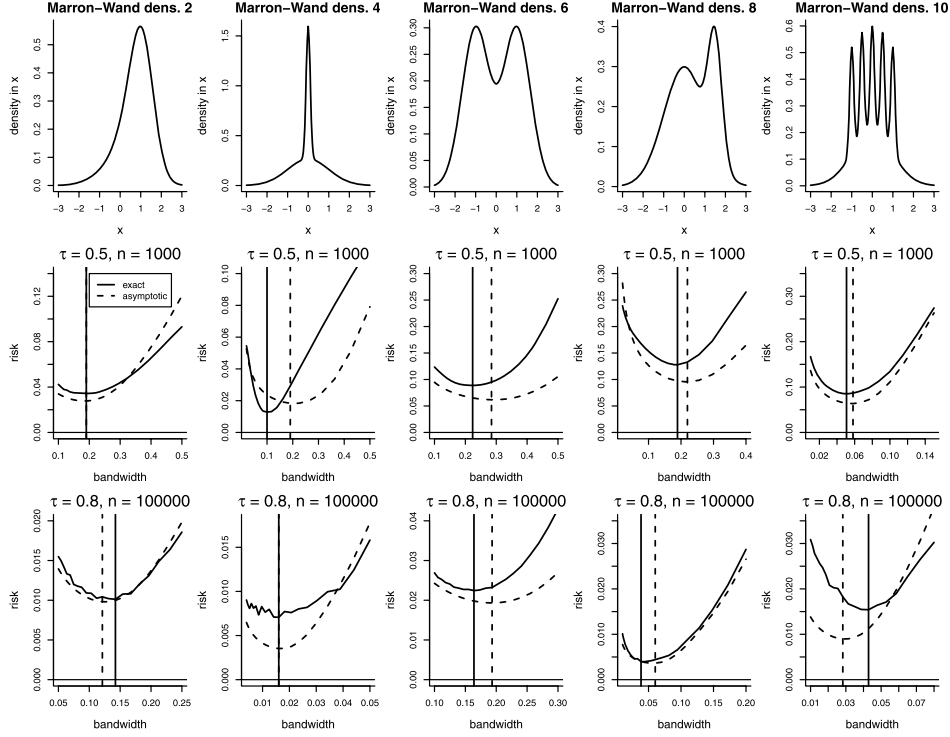}

\caption{Comparison of the ``exact'' risk
$\mathbb{E}\{\mu_f(\widehat{R}_{h,\tau} \triangle R_{\tau})\}$
and its asymptotic approximation for the five of the
Marron and Wand (\protect\citeyear{MW92}) density functions. The panels in
the first row are the density functions, and panels in the
same column correspond to the same density function.
In each panel in the second and third row, the ``exact'' risk,
obtained by averaging 100 realisations of
$\mu_f(\widehat{R}_{h,\tau} \triangle R_{\tau})$,
is shown as a solid black curve. The dashed curve is the asymptotic
risk approximation corresponding to the right-hand side of
(\protect\ref{eq:asyRiskSim}). Vertical lines pass through the minima
of the
``exact'' risk (solid line) and the asymptotic risk (broken line).
The second row has $\tau=0.5$ and $n=1000$, while the third row
has $\tau=0.8$ and $n=100\mbox{,}000$.}
\label{fig:HDRvAHDR}
\end{figure}

Figure \ref{fig:HDRvAHDR} compares the asymptotic risk approximation
with its ``exact'' counterpart for $f$ corresponding to Densities 2, 4,
6, 8 and 10 of Marron and Wand (\citeyear{MW92}), and $(\tau ,n) =
(0.5,1000)$ and $(0.8,100\mbox{,}000)$. The kernel $K$ is set to $\phi$
throughout, and the Monte Carlo sample size is $M=100$. For most of
these densities the asymptotic risk approximation is quite good for
$n=1000$ in the bandwidth range of interest. Density 4 is the main
exception; it appears that larger sample sizes are required for the
leading terms to be dominant. In particular, for this density, the
difficulty appears to be caused by very large values of $|f''|$ at the
crossing points of $f_{0.5}$ (for Density 4, the level $f_{0.5}$ is
very close to the rapid transition from shallow to steep gradient seen
in the corresponding upper panel in Figure \ref{fig:HDRvAHDR}). For
several densities, the estimand $R_{0.8}$ corresponds to the fine
detail of $f$. It is perhaps surprising that even with the larger
sample size, the asymptotic risk approximation is not always that
accurate, though in some cases the approximation is very good indeed.

\section{Bandwidth selection}\label{sec:bwSel}

In this section, we assume that, as in Corollary \ref{Cor:Opth}, there
exists a unique $c_{\mathrm{opt}} \in(0,\infty)$ such that any
optimal bandwidth sequence satisfies $h_{\mathrm{opt}} =
c_{\mathrm{opt}}n^{-1/5}\{1 + o(1)\}$. In this case,
$c_{\mathrm{opt}}$ minimizes the asymptotic risk given by
%
%
\begin{equation}
\label{Eq:copt}\quad
\mathrm{AR}(c) = \frac{1}{n^{2/5}}\sum_{j=1}^{2r} \biggl[\frac{B_{1,j}}{c^{1/2}}
\phi(B_{2,j}
c^{5/2}) + B_{3,j} c^2 \{2\Phi(B_{2,j} c^{5/2}) - 1\} \biggr].
\end{equation}
In order to find a practical bandwidth selector, we seek an estimator
$\widehat{c}_{\mathrm{opt}}$ of $c_{\mathrm{opt}}$. The natural way to
construct such an estimator is by using estimators $\widehat{D}_1$,
$\widehat{D}_2$ and $\widehat{D}_{3,j}$ of $D_1$, $D_2$ and $D_{3,j}$,
respectively, to obtain plug-in estimators $\widehat{B}_{1,j}$,
$\widehat{B}_{2,j}$ and $\widehat{B}_{3,j}$ of $B_{1,j}$, $B_{2,j}$ and
$B_{3,j}$, respectively. These\vspace*{1pt} in turn can be used to find $\widehat
{c}_{\mathrm{opt}}
= \arg\min_{c \in(0,\infty)} \widehat{\mathrm{AR}}_n(c)$, where
%
%
\begin{equation}
\label{Eq:hatcopt}\quad
\widehat{\mathrm{AR}}_n(c) = \frac{1}{n^{2/5}}\sum_{j=1}^{2r} \biggl[\frac{\widehat
{B}_{1,j}}{c^{1/2}} \phi(\widehat{B}_{2,j}
c^{5/2}) + \widehat{B}_{3,j} c^2 \{2\Phi(\widehat{B}_{2,j} c^{5/2}) - 1\} \biggr].
\end{equation}
With probability one, the solution to this minimization problem will
be unique for large $n$ provided that $\mathrm{AR}''(c_{\mathrm{opt}}) > 0$ and
this solution can easily be found numerically. Our final bandwidth
selector is then ${\widehat h}_{\tau \mathrm{HDR}}= \widehat{c}_{\mathrm
{opt}}n^{-1/5}$.

Note that we have not yet described how to construct the estimators
$\widehat{D}_1$, $\widehat{D}_2$ and~$\widehat{D}_{3,j}$. Again, we propose
plug-in estimators based on estimates of $f_{\tau }$ as well as
$f'(x_j)$ and $f''(x_j)$ for $j=1,\ldots,2r$. We assume the kernel
$K$ is smooth, and will construct kernel estimators
$\widehat{f}_{h_0}(\widehat{x}_{j,h_0})$, $\widehat{f}_{h_1}'(\widehat{x}_{j,h_0})$
and $\widehat{f}_{h_2}''(\widehat{x}_{j,h_0})$ of $f_{\tau }$, $f'(x_j)$
and $f''(x_j)$, respectively, where $\widehat{x}_{j,h_0}$ is an estimator
of $x_j$ described below. For the time being, we will use the same
kernel $K$ in all cases; this requirement will be dropped later on.
Even at this stage it will, however, be important to note that we can
use different bandwidths $h_0$, $h_1$ and $h_2$. Recall [e.g., Wand
and Jones (\citeyear{WJ95}), page 49] that,\vspace*{1pt} under appropriate conditions, if $h_k
\asymp n^{-1/(2k+5)}$ as $n \rightarrow\infty$ then
$\widehat{f}_{h_k}^{(k)}(x_j) - f^{(k)}(x_j) = O_p(n^{-2/(2k+5)})$ and
that this order cannot be improved for a nonnegative kernel. Here we
have used the notation $a_n \asymp b_n$ as $n \rightarrow\infty$ to
mean $0 < {\liminf_{n \rightarrow\infty}} |a_n/b_n| \leq{\limsup_{n
\rightarrow\infty}} |a_n/b_n| < \infty$. Finally, we observe that
if $h_0$ satisfies (A2), then with probability one, for all
sufficiently large $n$ there exist $\widehat{x}_{1,h_0} < \cdots<
\widehat{x}_{2r,h_0}$ such that $\widehat{f}_{h_0}(\widehat{x}_{j,h_0}) =
\widehat{f}_{h_0,\tau }$ for each $j$, and we use $\widehat{x}_{j,h_0}$ to
estimate $x_j$. Our theoretical study of the performance of this
bandwidth selector requires some additional conditions:
\begin{enumerate}[(A4):]
\item[(A4):] $f$ has four continuous derivatives in an open
set containing each $x_j$;
\item[(A5):] $h_0 \asymp n^{-1/5}$, $h_1 \asymp n^{-1/7}$ and
$h_2 \asymp n^{-1/9}$ as $n \rightarrow\infty$;
\item[(A6):] $K$ has a bounded third derivative, $K''$ is of
bounded variation and $\int|x|^3\times|K'(x)|+ x^4|K''(x)|\, dx <
\infty$.
\end{enumerate}
\begin{theorem}
\label{Thm:Hath}
Assume \textup{(A1)} and \textup{(A3)--(A6)}. Assume further
that $c_{opt}$ is unique and that $\mathrm{AR}''(c_{\mathrm{opt}}) > 0$. Then
\[
\frac{{\widehat h}_{\tau \mathrm{HDR}}}{h_{\mathrm{opt}}} = 1 +
O_p(n^{-2/9})
\]
as $n \rightarrow\infty$. Moreover, recalling that
${\widehat h}_{\tau \mathrm{HDR}}= \widehat{c}_{\mathrm{opt}}n^{-1/5}$, we
have
\[
\frac{\widehat{\mathrm{AR}}_n(\widehat{c}_{\mathrm{opt}})}{\mathrm{AR}(c_{\mathrm
{opt}})} = 1
+ O_p(n^{-2/9}).
\]
\end{theorem}

Examining the proof of Theorem \ref{Thm:Hath} reveals that the rate of
convergence to zero of the relative error of ${\widehat h}_{\tau
\mathrm{HDR}}$ is determined
by the rate at which we can estimate $f''(x_j)$ for $j=1,\ldots,2r$.
This suggests that we might be able to obtain a faster rate of
convergence by using a higher order kernel to estimate $f''(x_j)$ [and
in fact $f'(x_j)$]. Recall that we call $K$ an \textit{$S$th order
kernel} if:
\begin{enumerate}
\item$\int K(x)\, dx = 1$;
\item$\mu_s(K) \equiv\int x^s K(x) \,dx = 0$ for $s=1,\ldots,S-1$;
\item$\mu_S(K) \equiv\int x^S K(x) \,dx \neq0$ and $\int|x|^S |K(x)|\,
dx < \infty$.
\end{enumerate}
Higher order kernels refer to $S > 2$. The usual objection to the use
of higher order kernels, namely that such a kernel cannot be
nonnegative, is less significant when the aim is to estimate
derivatives of a density rather than the density itself. Let the
kernels used to estimate $f'(x_j)$ and $f''(x_j)$ be denoted $K_1$ and
$K_2$, respectively, and continue to denote the respective bandwidths
by $h_1$ and $h_2$. An improved rate of convergence of the relative
error of our bandwidth selector can be obtained by replacing
conditions (A4), (A5) and (A6) with the
following:
\begin{enumerate}[(A7):]
\item[(A7):] $f$ has 12 continuous derivatives in an open
set containing each $x_j$.
\item[(A8):] $h_0 \asymp n^{-1/5}$, $h_1 \asymp n^{-1/15}$ and
$h_2 \asymp n^{-1/25}$ as $n \rightarrow\infty$.
\item[(A9):] $K_1$ is a $6$th order kernel and has a bounded
second derivative with $K_1$ and $K_1'$ of bounded variation and
satisfying $\int x^6 |K_1(x)| + |x|^7|K_1'(x)| \,dx < \infty$.
Moreover, $K_2$ is a $10$th order kernel and has a bounded third
derivative with $K_2$, $K_2'$ and $K_2''$ of bounded variation and
satisfying $\int x^{10} |K_2(x)| + |x|^{11} |K_2'(x)| + x^{12}
|K_2''(x)| \,dx < \infty$.
\end{enumerate}
We write $\dhhat_{\tau \mathrm{HDR}}$ for the bandwidth selector
obtained in a similar
way to ${\widehat h}_{\tau \mathrm{HDR}}$, but using the kernels
$K_1$ and $K_2$ to estimate
$f'(x_j)$
and $f''(x_j)$, respectively, in the definitions of $D_1$, $D_2$, $D_{3,j}$,
$B_{1,j}$, $B_{2,j}$ and $B_{3,j}$.
\begin{theorem}
\label{Thm:HathHO}
Assume \textup{(A1)}, \textup{(A3)} and \textup{(A7)--(A9)}.
Assume further that $c_{opt}$ is unique and that
$\mathrm{AR}''(c_{\mathrm{opt}}) > 0$. Then
\[
\frac{\dhhat_{\tau \mathrm{HDR}}}{h_{\mathrm{opt}}} = 1 + O_p(n^{-2/5})
\]
as $n \rightarrow\infty$. Moreover, writing
$\dhhat_{\tau \mathrm{HDR}}= \widehat{\hspace*{-0.5pt}\widehat{c}}_{\mathrm{opt}}n^{-1/5}$,
we have
\[
\frac{\widehat{\mathrm{AR}}_n(\widehat{\hspace*{-0.5pt}\widehat{c}}_{\mathrm{opt}})}{\mathrm{AR}(c_{\mathrm
{opt}})} = 1 + O_p(n^{-2/5}).
\]
\end{theorem}

It is clear that Theorem \ref{Thm:Hath} represents a relatively weak
conclusion under relatively weak conditions, while Theorem \ref{Thm:HathHO}
represents a stronger conclusion under strong conditions. Intermediate
results are also possible but seem to be of little practical
interest.

\subsection{An effective practical bandwidth selector}\label{sec:pracAlg}

We confine our development of a practical consistent bandwidth
selector to the scenario where $f$ satisfies weaker smoothness
conditions of Theorem \ref{Thm:Hath}. Our end-product is a
fast-to-compute bandwidth selector for HDR estimation that possesses
the asymptotic properties conveyed by Theorem \ref{Thm:Hath}, performs
well in simulations and is readily implemented in~\texttt{R}. Indeed, as
detailed below, an \texttt{R} function for our procedure is available on
the Internet.

The pilot bandwidths $h_0$, $h_1$ and $h_2$ are estimated using
direct plug-in strategies with two levels of kernel functional
estimation. Chapter 3 of Wand and Jones (\citeyear{WJ95}) provides details on
this general approach to bandwidth selection. In the case of $h_0$ the
approach is similar to those proposed by Park and Marron (\citeyear{PM90}) and
Sheather and Jones (\citeyear{SJ91}). Direct plug-in bandwidth selection
strategies for density functions and their derivatives involve
estimation of functionals of the form
\[
\psi_r=\int_{-\infty}^{\infty} f^{(r)}(x)f(x) \,dx.
\]
Kernel estimators of $\psi_r$ take the form
\[
{\widehat\psi}_r(g)=n^{-2}g^{-r-1}\sum_{i=1}^n\sum_{j=1}^nL^{(r)}\{
(X_i-X_j)/g\},
\]
where $L$ is a sufficiently smooth 2nd-order kernel function, and $g>0$
is a
bandwidth parameter. Multi-level plug-in strategies use the fact
that the asymptotically optimal $g$, with respect\vspace*{1pt} to the mean
squared error of ${\widehat\psi}_r(g)$, is
$[-2L^{(r)}(0)/\{n\psi_{r+2}\int u^2L(u) \,du\}]^{1/(r+3)}$.
To get the algorithm started we also require
\textit{normal scale} estimates of $\psi_r$, based on the assumption
that $f$ is a $N(\mu,\sigma^2)$ density. Normal scale estimates
of $\psi_r$ take the form
\[
{\widehat\psi}^{\mathrm{NS}}_r=\frac{(-1)^{r/2}r!}{(2{\widehat
\sigma})^{r+1}(r/2)!\pi^{1/2}}.
\]
Throughout we take $K=L=\phi$, the standard normal kernel.
The full algorithm is:

\begin{enumerate}[10.]
\item The inputs are the random sample $X_1,\ldots,X_n$ and parameter
$0<\tau <1$ specifying the required HDR.
\item Let ${\widehat\sigma}=\min(\mbox{sample standard deviation},
(\mbox{sample interquartile range})/1.349)$ be
a robust estimate of scale. (The interquartile range for the standard
normal density is approximately $1.349$, so this factor ensures approximate
unbiasedness for normally distributed data.)
\item Estimate $\psi_8$, $\psi_{10}$ and $\psi_{12}$ using
normal scale estimates.\vspace*{1pt} Explicit expressions for these are
${\widehat\psi}^{\mathrm{NS}}_8=105/(32\pi^{1/2}{\widehat\sigma}^9)$,
${\widehat\psi}^{\mathrm{NS}}_{10}=-945/(64\pi^{1/2}{\widehat
\sigma}^{11})$
and
${\widehat\psi}^{\mathrm{NS}}_{12}=10395/(128\pi^{1/2}{\widehat
\sigma}^{13})$.

\item Estimate\vspace*{1pt} $\psi_6$, $\psi_8$ and $\psi_{10}$ using
kernel estimates ${\widehat\psi}_6(g_{0,1})$,
${\widehat\psi}_8(g_{1,1})$ and ${\widehat\psi}_{10}(g_{1,1})$ where
$g_{0,1}=\{30/({\widehat\psi}^{\mathrm{NS}}_8 n)\}^{1/9}$,
$g_{1,1}=\{-210/({\widehat\psi}^{\mathrm{NS}}_{10} n)\}^{1/11}$
and
$g_{1,2}=\{1890/({\widehat\psi}^{\mathrm{NS}}_{12} n)\}^{1/13}$.
\item
Estimate\vspace*{1pt} $\psi_4$, $\psi_6$ and $\psi_8$ using
kernel estimates ${\widehat\psi}_4(g_{0,2})$,
${\widehat\psi}_6(g_{1,2})$ and ${\widehat\psi}_8(g_{2,2})$
where\vspace*{1pt}
$g_{0,2}=[6/\{{\widehat\psi}_8(g_{0,1}) n\}]^{1/7}$,
$g_{1,2}=[-30/\{{\widehat\psi}_{10}(g_{1,1}) n\}]^{1/9}$
and
$g_{2,2}=[210/\{{\widehat\psi}_{12}(g_{1,2}) n\}]^{1/11}$.
\item Obtain direct plug-in bandwidths ${\widehat h}^{(r)}$
for estimation of $f^{(r)}$ by replacing $\psi_{r+2}$
in the optimal expression, with respect\vspace*{1pt}
to asymptotic mean integrated squared error,
by ${\widehat\psi}_{r+2}(g_{r,2})$. Explicit expressions\vspace*{1pt} are
${\widehat h}_0=[1/\{2\pi^{1/2}{\widehat\psi}_4(g_{0,2})n\}]^{1/5}$,
${\widehat h}_1=[-3/\{4\pi^{1/2}{\widehat\psi}_6(g_{1,2})n\}]^{1/7}$
and ${\widehat h}_2=[15/\{8\pi^{1/2}{\widehat\psi}_8(g_{2,2})n\}]^{1/9}$.
\item Obtain pilot of estimates of $f$, $f'$ and $f''$ via Gaussian
kernel estimates based on these bandwidths: ${\widehat f}_{{\widehat
h}_0}(\cdot)$,
${\widehat f}'_{{\widehat h}_1}(\cdot)$ and ${\widehat f}_{{\widehat
h\,}_2}''(\cdot)$.
\item
Use ${\widehat f}_{{\widehat h}_0}(\cdot)$ to obtain pilot estimates of
$f_{\tau }$, $r$ and $x_1,\ldots,x_{2r}$.
\item Substitute the estimates from Steps \ref{Step6} and \ref{Step7} into
the expressions for $B_{1,j}$, $B_{2,j}$ and $B_{3,j}$ to
obtain estimates ${\widehat B}_{1,j}$, ${\widehat B}_{2,j}$ and
${\widehat B}_{3,j}$.
\item The selected bandwidth for Gaussian kernel estimation of the
$100(1-\tau )\%$ HDR is ${\widehat h}_{\tau \mathrm{HDR}}=\widehat
{c}_{\mathrm{opt}} n^{-1/5}$
where $\widehat{c}_{\mathrm{opt}} = \arg\min_{c \in(0,\infty)}
\widehat{\mathrm{AR}}_n(c)$, where $\widehat{\mathrm{AR}}_n$ was defined in (\ref{Eq:hatcopt}).
\end{enumerate}

Binned approximations to ${\widehat\psi}_r(g)$ [cf. Gonz\'alez-Manteiga,
Sanch\'ez-Sellero and Wand (\citeyear{GSW96})] are strongly recommended to allow
fast processing of large samples. An \texttt{R} function \texttt{hdrbw()}
that implements the above algorithm has been included in the package
\texttt{hdrcde} [Hyndman (\citeyear{H09})] which supports HDR estimation.

\subsection{Simulation results}

\begin{figure}

\includegraphics{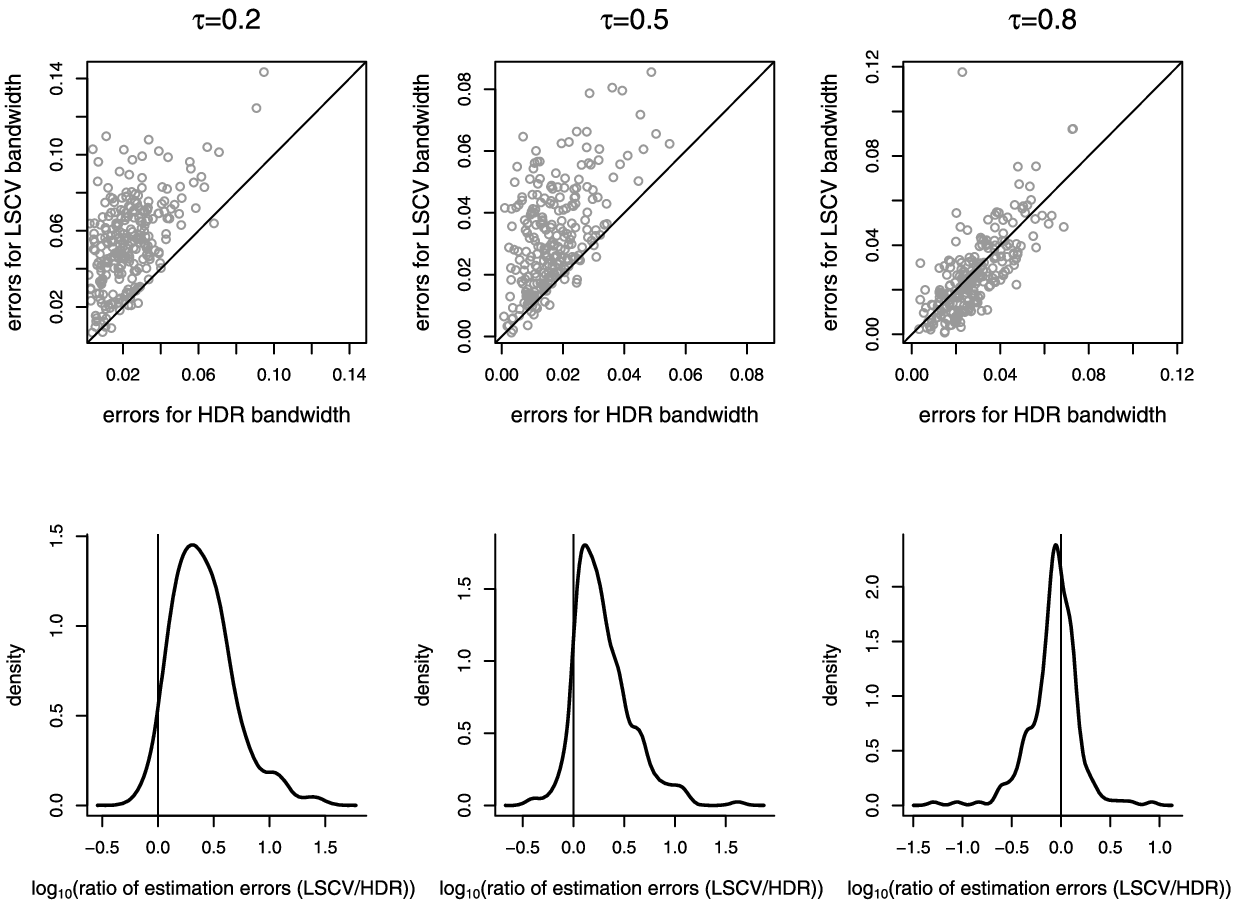}

\caption{Summary of simulation comparison between
${\widehat h}_{\tau \mathrm{HDR}}$ and ${\widehat h}_{\mathrm
{LSCV}}$ for $\tau =0.2,0.5$ and $0.8$ and
250 samples of size 1000 generated from Density 4 of
Marron and Wand (\protect\citeyear{MW92}).
The upper panels are scatterplots of the errors
$\mu_f(\widehat{R}_{h,\tau } \triangle R_\tau )$
for $h={\widehat h}_{\mathrm{LSCV}}$ on the vertical axes and
$h={\widehat h}_{\tau \mathrm{HDR}}$
on the horizontal axes.
The lower panels are kernel density estimates of
$\log_{10}((\mbox{error for } h={\widehat h}_{\tau \mathrm
{HDR}})/(\mbox{error for }h={\widehat h}_{\mathrm{LSCV}}))$.}
\label{fig:bwSimD4n1000}
\end{figure}

\begin{figure}

\includegraphics{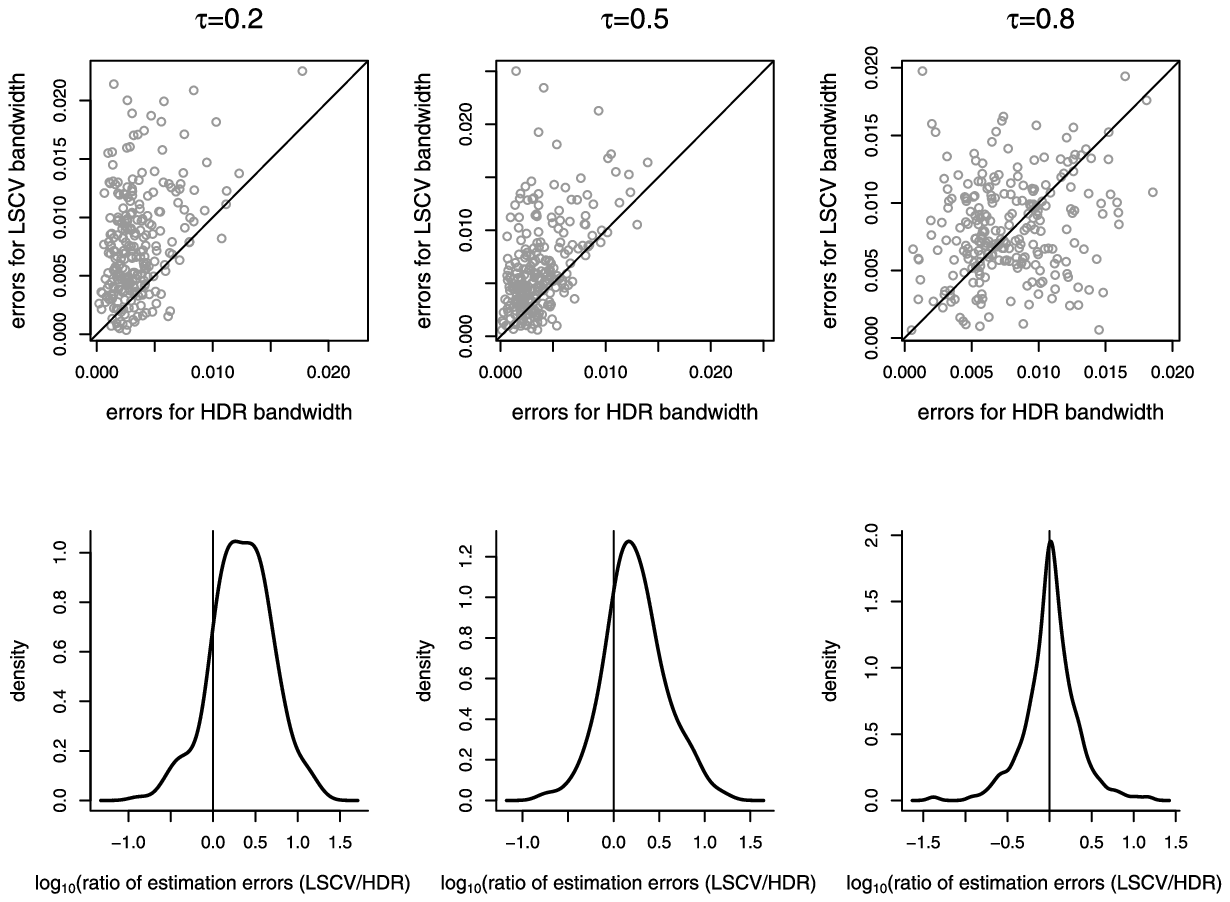}

\caption{Summary of simulation comparison between
${\widehat h}_{\tau \mathrm{HDR}}$ and ${\widehat h}_{\mathrm
{LSCV}}$ for $\tau =0.2,0.5$ and $0.8$
and 250 samples of size 100\mbox{,}000 generated from Density 4
of Marron and Wand (\protect\citeyear{MW92}). The upper panels
are scatterplots of the errors
$\mu_f(\widehat{R}_{h,\tau } \triangle R_\tau )$
for $h={\widehat h}_{\mathrm{LSCV}}$ on the vertical axes and
$h={\widehat h}_{\tau \mathrm{HDR}}$
on the horizontal axes.
The lower panels are kernel density estimates of
$\log_{10}((\mbox{error for }h={\widehat h}_{\tau \mathrm
{HDR}})/(\mbox{error for }h={\widehat h}_{\mathrm{LSCV}}))$.}
\label{fig:bwSimD4n100000}
\end{figure}

We ran a simulation study in which the performance of ${\widehat
h}_{\tau \mathrm{HDR}}$
was compared with an established ISE-based selector: least
squares\vspace*{1pt}
cross validation [Rudemo (\citeyear{R82}), Bowman (\citeyear{B84})] which we denote by
${\widehat h}_{\mathrm{LSCV}}$. The number of replications in the
simulation study was
250. The HDR estimation error $\mu_f(\widehat{R}_{h,\tau } \triangle
R_\tau )$ was used throughout the study. Figures
\ref{fig:bwSimD4n1000} ($n=1000$) and \ref{fig:bwSimD4n100000}
($n=100\mbox{,}000$) summarise the results for the situation where the true
$f$ is the normal mixture density from Section \ref{sec:intro} and
Figure \ref{fig:HDRvISE}. The improvement gained from using the
HDR-tailored bandwidth selector is apparent from the graphics,
especially for the lower values of $\tau $. Wilcoxon tests applied
to the error ratios showed statistically significant improvement of
${\widehat h}_{\tau \mathrm{HDR}}$ at the 5\% level for $\tau
=0.2,0.5$ and $0.8$ when
$n=100\mbox{,}000$. For $n=1000$, ${\widehat h}_{\tau \mathrm{HDR}}$
performed better for
$\tau =0.2,0.5$, while ${\widehat h}_{\mathrm{LSCV}}$ did better for
$\tau =0.8$.
This latter result is not a big surprise since good estimation of
$R_{0.8}$ requires good estimation of the finger-shaped modal region
and this, in turn, requires good ISE performance.

%

%

We performed similar simulation comparisons for the
remaining Densities 1--10 of Marron and Wand (\citeyear{MW92}).
For $n=1000$ the performance of ${\widehat h}_{\tau \mathrm{HDR}}$
was better
than ${\widehat h}_{\mathrm{LSCV}}$ for Densities 1--5; whereas
${\widehat h}_{\mathrm{LSCV}}$
did better for Densities 6--10. This suggests that
the asymptotics on which ${\widehat h}_{\tau \mathrm{HDR}}$ relies
have not
``kicked in'' at $n=1000$ for these more intricate density functions.
We suspect that more sophisticated pilot estimation might
improve matters for HDR-based bandwidth selection for
lower sample sizes.
The $n=100\mbox{,}000$ simulations show superior performance
of ${\widehat h}_{\tau \mathrm{HDR}}$, especially $\tau =0.8$ where it
is the ``winner'' for 9 out of the 10 density functions.
The overarching conclusion is that for common
density estimation situations ${\widehat h}_{\tau \mathrm{HDR}}$ is better
than ${\widehat h}_{\mathrm{LSCV}}$.

\subsection{Application to daily temperature data}

We conclude with an application to data on daily maximum temperatures
in Melbourne, Australia, for the years 1981--1990. These data were used
in Hyndman (\citeyear{H96}) to illustrate HDR principles. We revisit them armed with
the automatic HDR estimation technology described in Section \ref{sec:pracAlg}.
Of interest are the conditional densities of
tomorrow's temperature \textit{given} today's temperature
is within a fixed interval.

\begin{figure}

\includegraphics{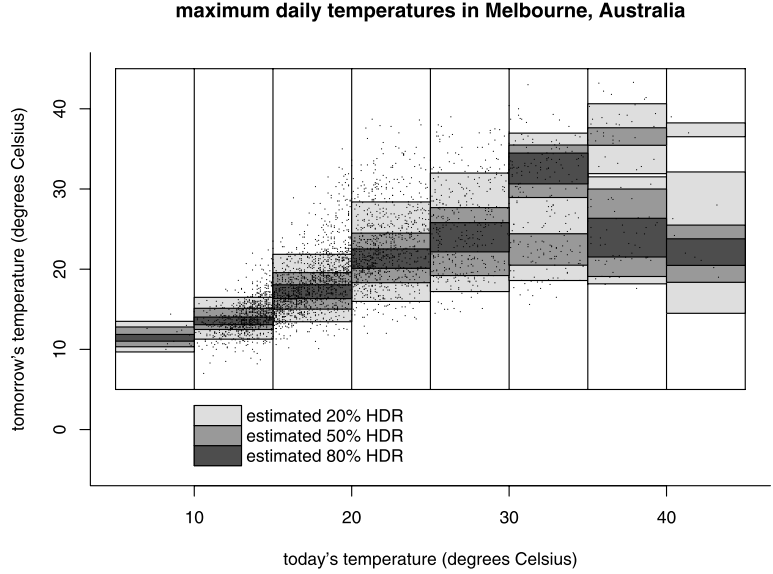}

\caption{Estimated kernel HDRs for the conditional densities
of tomorrow's temperature given that today's temperature
is in a fixed interval. The bandwidth for each
HDR estimate is chosen using the selector described in Section \protect
\ref{sec:pracAlg}.}
\label{fig:melbmaxHDR}
\end{figure}

The intervals for the ``today's temperature'' values are, in degrees Celsius,
\[
[5,10),[10,15),\ldots,[40,45).
\]
Figure \ref{fig:melbmaxHDR} shows the kernel
%
%
estimates of the 20\%, 50\% and 80\% HDRs with bandwidths chosen using
the rule ${\widehat h}_{\tau \mathrm{HDR}}$ as detailed in Section
\ref{sec:pracAlg}.
Some interesting bimodality in ``tomorrow's temperature'' is apparent when
conditioned on today's temperature being in the 30--40 degrees Celsius range.

\begin{appendix}\label{app}
\section*{Appendix: Proofs}

\subsection*{Proof of Theorem \protect\ref{Thm:Mainthm}}

Throughout the proof, it is convenient to write $x_0 = -\infty$ and
$x_{2r+1} = \infty$
and adopt the convention that $x_0 + a = -\infty$ and $x_{2r+1} + a =
\infty$ for
all $a \in\mathbb{R}$. Observe that
\begin{eqnarray*}
\mu_f(\widehat{R}_{h,\tau } \triangle R_\tau )
&=& \int_{-\infty}^\infty f(x) \bigl|
\mathbbm{1}_{\{\widehat{f}_h(x) \geq\widehat{f}_{h,\tau }\}}
- \mathbbm{1}_{\{f(x) \geq f_\tau \}}\bigr| \,dx \\
&=& \sum_{j=0}^r \int_{x_{2j}}^{x_{2j+1}} f(x)
\mathbbm{1}_{\{\widehat{f}_h(x) \geq\widehat{f}_{h,\tau }\}} \,dx
\\
&&{}
+ \sum_{j=1}^r \int_{x_{2j-1}}^{x_{2j}} f(x)
\mathbbm{1}_{\{\widehat{f}_h(x) < \widehat{f}_{h,\tau }\}} \,dx,
\end{eqnarray*}
so that
\begin{eqnarray*}
\mathbb{E}\{\mu_f(\widehat{R}_{h,\tau } \triangle R_\tau )\} &=&
\sum_{j=0}^r \int_{x_{2j}}^{x_{2j+1}} f(x)
\mathbb{P}\bigl(\widehat{f}_h(x) \geq\widehat{f}_{h,\tau }\bigr) \,dx\\
&&{} +
\sum_{j=1}^r \int_{x_{2j-1}}^{x_{2j}} f(x)
\mathbb{P}\bigl(\widehat{f}_h(x) < \widehat{f}_{h,\tau }\bigr) \,dx.
\end{eqnarray*}
The main idea of the proof is that the dominant contribution to
$\mathbb{E}\{\mu_f(\widehat{R}_{h,\tau } \triangle R_\tau )\}$ comes
from a union of $2r$ small intervals, one near each $x_j$, where
$\mathbb{P}(\widehat{f}_h(x) \geq\widehat{f}_{h,\tau })$ is close
to $1/2$. In each of these intervals, we can represent $\widehat{f}_h(x)
- \widehat{f}_{h,\tau }$ by a sample mean of independent and identically
distributed random variables and a small additional remainder term,
and hence apply a normal approximation to deduce the result. For
clarity of exposition, we now split the proof into several steps:
\begin{Step}\label{Step1}
As a preliminary step, let $\tilde{f} = f+g$ be
another uniformly continuous density, and let $\tilde{f}_{\tau } =
f_{\tau }(\tilde{f})$. Writing $\|\cdot\|_{\infty}$ for the supremum
norm on the real line, we show that there exists $C \geq1$ such that
for all $\varepsilon> 0$ sufficiently small, we have
$|\tilde{f}_{\tau } - f_{\tau }| \leq C \varepsilon$ whenever
$\|g\|_{\infty} \equiv\|\tilde{f} - f\|_\infty\leq\varepsilon$. To
see this, let $L = \sum_{j=1}^r (x_{2j} - x_{2j-1})$ and choose $C >
1+ 2L\{\frac{1}{4}f_{\tau } \sum_{j=1}^{2r}
\frac{1}{|f'(x_j)|}\}^{-1}$. The inverse function theorem
[Burkill and Burkill (\citeyear{BB02}), Theorem 7.51] gives that for $\varepsilon\in
\mathbb{R}$ with $|\varepsilon|$ sufficiently small, we can write
\[
\{x\dvtx f(x) \geq f_{\tau } + \varepsilon\} = \bigcup_{j=1}^r
[x_{2j-1} +
\delta_{\varepsilon,2j-1},x_{2j} - \delta_{\varepsilon,2j}]
\]
with $\delta_{\varepsilon,j} = \frac{\varepsilon}{|f'(x_j)|} +
O(\varepsilon^2)$ as $\varepsilon\rightarrow0$. It follows that when
$\varepsilon> 0$ is sufficiently small, and $\|g\|_\infty\leq
\varepsilon$, we have
\begin{eqnarray*}
&&\int_{-\infty}^\infty\tilde{f}(x) \mathbbm{1}_{\{\tilde{f}(x)
\geq
f_\tau - C\varepsilon\}} \,dx\\
&&\qquad\geq\int_{-\infty}^\infty\{f(x) - \varepsilon\} \mathbbm{1}_{\{f(x)
\geq
f_\tau
- (C-1)\varepsilon\}} \,dx \\
&&\qquad= 1-\tau + \sum_{j=1}^r \int_{x_{2j-1} +
\delta_{-\varepsilon(C-1),2j-1}}^{x_{2j-1}} f(x) \,dx\\
&&\qquad\quad{} + \sum_{j=1}^r
\int_{x_{2j}}^{x_{2j} - \delta_{-\varepsilon(C-1),2j}} f(x) \,dx \\
&&\qquad\quad{} - \varepsilon\sum_{j=1}^r \bigl\{x_{2j} - \delta_{-\varepsilon(C-1),2j}
- \bigl(x_{2j-1} + \delta_{-\varepsilon(C-1),2j-1}\bigr)\bigr\} \\
&&\qquad\geq 1-\tau + \frac{1}{4}(C-1)\varepsilon f_{\tau }\sum_{j=1}^{2r}
\frac{1}{|f'(x_j)|} - 2 \varepsilon L > 1-\tau .
\end{eqnarray*}
Thus $\tilde{f}_\tau \geq f_\tau - C \varepsilon$. A very similar
argument yields the upper bound $\tilde{f}_\tau \leq f_\tau
+ C
\varepsilon$, and this completes Step \ref{Step1}.
\end{Step}
\begin{Remark*} Now, for $\delta> 0$ small enough that $f$ has two
continuous derivatives in $I_\delta\equiv\bigcup_{j=1}^r
[x_{2j-1}-\delta,x_{2j}+\delta]$, let $\|\cdot\|_{I_\delta,\infty}$
denote the supremum norm restricted to $I_\delta$. It will be helpful
in Step \ref{Step4} to note that a small modification of the above
argument may be used to prove that if $\|g\|_\infty$ and
$\|g'\|_{I_\delta,\infty}$ are sufficiently small, and if
\[
\sum_{j=1}^r \int_{x_{2j-1}-\delta}^{x_{2j}+\delta} |g(x)| \,dx =
O \Biggl(\sum_{j=1}^{2r} |g(x_j)| \Biggr)
\]
as $\sum_{j=1}^{2r} |g(x_j)| \rightarrow0$, then $\tilde
{f}_\tau -
f_\tau = O (\sum_{j=1}^{2r} |g(x_j)| )$ as $\sum_{j=1}^{2r}
|g(x_j)| \rightarrow0$.
\end{Remark*}
\begin{Step}\label{Step2}
We show that for each fixed $\delta> 0$,
%
%
\setcounter{equation}{0}
\begin{eqnarray}
\label{Eq:NonMargin}
&&\sum_{j=0}^r \int_{x_{2j}+\delta}^{x_{2j+1}-\delta} f(x) \mathbb
{P}\bigl(\widehat{f}_h(x)
\geq\widehat{f}_{h,\tau }\bigr) \,dx\nonumber\\[-8pt]\\[-8pt]
&&\qquad{} + \sum_{j=1}^r \int_{x_{2j-1}+\delta
}^{x_{2j}-\delta}
f(x) \mathbb{P}\bigl(\widehat{f}_h(x) < \widehat{f}_{h,\tau }\bigr) \,dx =
o(n^{-1})\nonumber
\end{eqnarray}
as $n \rightarrow\infty$. In fact, we claim (and it will be
straightforward to see) that the error term is of the stated order
uniformly for $h \in[h^-,h^+]$. Indeed, we make a similar claim for
every error term in each expression below where the bandwidth $h$
appears, but we do not repeat this assertion in future occurrences.
As in Step \ref{Step1}, observe that under (A1), if $\delta>
0$ is sufficiently small, then there exists $\varepsilon> 0$ such that
$f(x) \leq f_\tau - \varepsilon$ for $x \in\bigcup_{j=0}^r
[x_{2j}+\delta,x_{2j+1}-\delta]$ and $f(x) \geq f_\tau +
\varepsilon$
for $x \in\bigcup_{j=1}^r [x_{2j-1}+\delta,x_{2j}-\delta]$. By reducing
$\delta> 0$ if necessary, for $x \in\bigcup_{j=0}^r
[x_{2j}+\delta,x_{2j+1}-\delta]$,
%
%
\begin{eqnarray}
\label{Eq:BasicIneq}
\mathbb{P}\bigl(\widehat{f}_h(x) \geq\widehat{f}_{h,\tau }\bigr)
&\leq&\mathbb{P}\bigl(\widehat{f}_h(x) - f(x) - (\widehat{f}_{h,\tau } -
f_{\tau })
\geq\varepsilon\bigr) \nonumber\\
&\leq&\mathbb{P}(\|\widehat{f}_h - f\|_{\infty} \geq\varepsilon/2)
+ \mathbb{P}(|\widehat{f}_{h,\tau } - f_{\tau }| \geq\varepsilon/2)
\\
&\leq& 2\mathbb{P} \biggl(\|\widehat{f}_h - f\|_{\infty} \geq\frac{\varepsilon
}{2C} \biggr),\nonumber
\end{eqnarray}
where we have used the result of Step \ref{Step1} in the last
inequality, and $C$ is the constant defined in that step. A very
similar argument yields the same upper bound for
$\mathbb{P}(\widehat{f}_h(x) < \widehat{f}_{h,\tau })$ when $x \in
\bigcup_{j=1}^r [x_{2j-1}+\delta,x_{2j}-\delta]$. Now, since $f$ is
uniformly continuous under (A1),
%
%
\begin{equation}
\label{Eq:Bias}
\|\mathbb{E}(\widehat{f}_h) - f\|_{\infty} = \sup_{x \in\mathbb{R}}
\biggl|\int_{-\infty}^\infty K(z)\{f(x-hz) - f(x)\}\, dz \biggr| \rightarrow0
\end{equation}
as $n \rightarrow\infty$. The inequality (\ref{Eq:BasicIneq}),
together with the observation (\ref{Eq:Bias}) on the bias of
$\widehat{f}_h$, yields that for $n$ sufficiently large,
\begin{eqnarray*}
&&\sum_{j=0}^r \int_{x_{2j}+\delta}^{x_{2j+1}-\delta} f(x)
\mathbb{P}\bigl(\widehat{f}_h(x) \geq\widehat{f}_{h,\tau }\bigr) \,dx\\
&&\quad{} + \sum_{j=1}^r \int_{x_{2j-1}+\delta}^{x_{2j}-\delta} f(x)
\mathbb{P}\bigl(\widehat{f}_h(x) < \widehat{f}_{h,\tau }\bigr) \,dx \\
&&\qquad\leq2\mathbb{P} \biggl(\|\widehat{f}_h -
\mathbb{E}(\widehat{f}_h)\|_{\infty} \geq\frac{\varepsilon}{4C}
\biggr)\\
&&\qquad\leq\exp(-c_1 nh\varepsilon^2)
\end{eqnarray*}
for some $c_1 > 0$. Here, the final inequality is an application of
Corollary 2.2 of Gin\'e and Guillou (\citeyear{GG02}) (a consequence of Talagrand's
inequality) to the Vapnik--Cervonenkis class of functions $\{K((x -
\cdot)/h)\dvtx x \in\mathbb{R}, h > 0\}$ [cf. Dudley (\citeyear{D99}),
Theorems 4.2.1 and 4.2.4]. Equation (\ref{Eq:NonMargin}) follows
immediately, and this completes the proof of Step \ref{Step2}.
\end{Step}
\begin{Step}\label{Step3}
We show that (\ref{Eq:NonMargin}) continues to hold
if $\delta$ is replaced by a sequence $(\delta_n)$ converging to zero,
provided that $\delta_n \rightarrow0$ slowly enough that
$n^{1/4}\delta_n \rightarrow\infty$ and $(h^+)^2 =
o(\delta_n)$. In order to complete the proof of Step \ref{Step3}, it
suffices to show that there exists $\delta> 0$ such that
\begin{eqnarray*}
E(\delta,\delta_n) &\equiv& \sum_{j=1}^r \int_{x_{2j-1} -
\delta}^{x_{2j-1} - \delta_n} \mathbb{P}\bigl(\widehat{f}_h(x) \geq
\widehat{f}_{h,\tau }\bigr) \,dx\\
&&{} + \sum_{j=1}^r \int_{x_{2j-1}
 +
\delta_n}^{x_{2j-1} + \delta} \mathbb{P}\bigl(\widehat{f}_h(x) <
\widehat{f}_{h,\tau }\bigr) \,dx
\sum_{j=1}^r
\int_{x_{2j} - \delta}^{x_{2j} - \delta_n} \mathbb{P}\bigl(\widehat{f}_h(x) <
\widehat{f}_{h,\tau }\bigr) \,dx\\
&&{} + \sum_{j=1}^r \int_{x_{2j}+
\delta_n}^{x_{2j} + \delta} \mathbb{P}\bigl(\widehat{f}_h(x) \geq
\widehat{f}_{h,\tau }\bigr) \,dx = o(n^{-1}).
\end{eqnarray*}
We may assume $\delta> 0$ is small enough that $f$ has two continuous
derivatives in~$I_{\delta}$. This enables a straightforward
modification to the argument in (\ref{Eq:Bias}) using a Taylor
expansion, leading to
%
%
\begin{equation}
\label{Eq:Bias2}
\|\mathbb{E}(\widehat{f}_h) - f\|_{I_\delta,\infty} = O(h^2).
\end{equation}
Now there exists a constant $c_2 > 0$ small enough that if we take
$\varepsilon_n = c_2\delta_n$, then we have $|f(x) - f_{\tau }|
\geq
\varepsilon_n$ when $\min_j |x-x_j| \geq\delta_n$. Moreover, $(h^+)^2 =
o(\varepsilon_n)$, so that for $n$ sufficiently large, the same argument
as in Step \ref{Step2} yields
\begin{eqnarray*}
E(\delta,\delta_n) &\leq& 2\mathbb{P} \biggl(\|\widehat{f}_h - \mathbb{E}(\widehat
{f}_h)\|_{\infty}
\geq\frac{\varepsilon_n}{4C} \biggr)\\ &\leq& \exp(-c_1 nh\varepsilon_n^2) =
o(n^{-1}).
\end{eqnarray*}
This completes the proof of Step \ref{Step3}.
\end{Step}
\begin{Step}\label{Step4}
We seek asymptotic expansions for
$\mathbb{E}(\widehat{f}_{h,\tau })$ and
$\operatorname{Var}(\widehat{f}_{h,\tau })$. To this end, for uniformly
continuous densities $\tilde{f} = f+g$ that are twice continuously
differentiable in $I_{\delta}$, and for $y
\in(0,\infty)$, we define
\[
\psi(\tilde{f},y) = \int_{-\infty}^\infty
\tilde{f}(x)\mathbbm{1}_{\{\tilde{f}(x) \geq y\}} \,dx.
\]
The reason for making this definition is that by examining the
behavior of $\psi$ under small changes of its arguments from
$(f,f_{\tau })$, we will be able to study the difference
$\widehat{f}_{h,\tau } - f_{\tau }$ in (\ref{Eq:ProbBound}) below.
First, for $\varepsilon> 0$ sufficiently small,
%
%
\begin{eqnarray}
\label{Eq:Partial1}
&&\Biggl|\psi(f,f_{\tau } +\varepsilon) - \psi(f,f_{\tau }) +
\varepsilon f_\tau \sum_{j=1}^{2r}\frac{1}{|f'(x_j)|} \Biggr| \nonumber\\
&&\qquad=
\Biggl|- \int_{-\infty}^\infty f(x) \mathbbm{1}_{\{f_{\tau } \leq
f(x) < f_{\tau }+\varepsilon\}} \,dx +
\varepsilon f_\tau \sum_{j=1}^{2r}\frac{1}{|f'(x_j)|} \Biggr|\nonumber\\
&&\qquad= \Biggl| -\sum_{j=1}^r \biggl\{\int_{x_{2j-1}}^{x_{2j-1}
+\delta_{\varepsilon,2j-1}} f(x) \,dx\\
&&\hspace*{75.7pt}{} +
\int_{x_{2j}-\delta_{\varepsilon,2j}}^{x_{2j}}
f(x) \,dx \biggr\} + \varepsilon f_\tau \sum_{j=1}^{2r}\frac{1}{|f'(x_j)|}
\Biggr|\nonumber\\
&&\qquad=O(\varepsilon^2)\nonumber
\end{eqnarray}
as $\varepsilon\searrow0$. A very similar argument shows that the
error term is of the same order as $\varepsilon\nearrow0$.

Observe that when $\|g\|_{\infty}$ and $\|g'\|_{I_{\delta},\infty}$
are sufficiently small, $\tilde{f}$ has a nonzero derivative in a
neighborhood of each $x_j$. It follows that for sufficiently small
values of $\|g\|_{\infty} + \|g'\|_{I_{\delta},\infty}$, we can write
\[
\{x\dvtx\tilde{f}(x) \geq\tilde{f}_{\tau }\} = \bigcup_{j=1}^r
[x_{2j-1} +
\delta_{\varepsilon,2j-1} + \eta_{2j-1} , x_{2j} -
\delta_{\varepsilon,2j} - \eta_{2j}],
\]
where $\varepsilon= \tilde{f}_\tau - f_\tau $. Moreover, provided
that
\[
\sum_{j=1}^r \int_{x_{2j-1}-\delta}^{x_{2j}+\delta} |g(x)| \,dx
= O\Biggl(\sum_{j=1}^{2r} |g(x_j)|\Biggr)
\]
and $\sum_{j=1}^{2r} |g(x_j)|
= O(\min_j |g(x_j)|)$ as $\sum_{j=1}^{2r} |g(x_j)| +
\|g'\|_{I_\delta,\infty} \rightarrow0$, we have that $\eta_j =
\frac{-g(x_j)}{|f'(x_j)|} + O(|g(x_j)|\|g'\|_{I_\delta,\infty})$ as
$\sum_{j=1}^{2r} |g(x_j)| + \|g'\|_{I_\delta,\infty} \rightarrow0$.
Thus we can write
%
%
\begin{eqnarray}
\label{Eq:Partial2}\quad
&&\Biggl|\psi(\tilde{f},\tilde{f}_\tau ) - \psi(f,\tilde
{f}_\tau ) -
f_{\tau }\sum_{j=1}^{2r} \frac{g(x_j)}{|f'(x_j)|} - \sum_{j=1}^r
\int_{x_{2j-1}}^{x_{2j}} g(x) \,dx \Biggr| \nonumber\\
&&\qquad\leq\Biggl|\int_{-\infty}^\infty f(x)\bigl(\mathbbm{1}_{\{\tilde{f}(x)
\geq\tilde{f}_\tau \}} - \mathbbm{1}_{\{f(x) \geq
\tilde{f}_\tau \}}\bigr) \,dx - f_{\tau }\sum_{j=1}^{2r}
\frac{g(x_j)}{|f'(x_j)|} \Biggr| \nonumber\\
&&\qquad\quad{} + \biggl|\int_{-\infty}^\infty g(x)
\bigl(\mathbbm{1}_{\{\tilde{f}(x) \geq\tilde{f}_\tau \}} -
\mathbbm{1}_{\{f(x) \geq f_\tau \}}\bigr) \,dx \biggr| \nonumber\\
&&\qquad= \Biggl| \Biggl\{f_{\tau } + O \Biggl(\sum_{j=1}^{2r} |g(x_j)| \Biggr) \Biggr\} \sum_{j=1}^{2r}
\Biggl\{\frac{g(x_j)}{|f'(x_j)|} +
O(|g(x_j)|\|g'\|_{I_\delta,\infty}) \Biggr\}\\
&&\hspace*{233pt}{} - f_{\tau }
\sum_{j=1}^{2r} \frac{g(x_j)}{|f'(x_j)|} \Biggr| \nonumber\\
&&\qquad\quad{} + O \Biggl\{ \Biggl(\sum_{j=1}^{2r} |g(x_j)| \Biggr)^2 \Biggr\}\nonumber\\
&&\qquad= O \Biggl\{ \Biggl(\sum_{j=1}^{2r} |g(x_j)| \Biggr)^2 + \|g'\|_{I_\delta,\infty}\sum
_{j=1}^{2r} |g(x_j)| \Biggr\}\nonumber
\end{eqnarray}
as $\sum_{j=1}^{2r} |g(x_j)| + \|g'\|_{I_\delta,\infty} \rightarrow
0$. Assuming that $\psi(\tilde{f},\tilde{f}_{\tau }) =
1-\tau $ and that the above conditions on $g$ hold, we have
from (\ref{Eq:Partial1}) and (\ref{Eq:Partial2}) that
%
%
\begin{eqnarray}
\label{Eq:PutTogether}\qquad
0 &=& \psi(\tilde{f},\tilde{f}_{\tau })
- \psi(f,f_{\tau }) \nonumber\\
&=& \psi(\tilde{f},\tilde{f}_{\tau }) -
\psi(f,\tilde{f}_{\tau })
+ \psi(f,\tilde{f}_{\tau }) - \psi(f,f_{\tau }) \nonumber
\nonumber\\[-8pt]\\[-8pt]
&=& - \{\tilde{f}_{\tau } - f_{\tau }\} f_{\tau }\sum_{j=1}^{2r}
\frac{1}{|f'(x_j)|} + f_{\tau }\sum_{j=1}^{2r}
\frac{g(x_j)}{|f'(x_j)|} + \sum_{j=1}^r \int_{x_{2j-1}}^{x_{2j}} g(x)
\,dx \nonumber\\
&&{} + O \Biggl\{ \Biggl(\sum_{j=1}^{2r} |g(x_j)| \Biggr)^2 +
\|g'\|_{I_\delta,\infty}\sum_{j=1}^{2r} |g(x_j)| \Biggr\}\nonumber
\end{eqnarray}
as $\sum_{j=1}^{2r} |g(x_j)| + \|g'\|_{I_\delta,\infty} \rightarrow
0$.

We want to apply (\ref{Eq:PutTogether}) with $\tilde{f} = \widehat{f}_h$,
so that $g = \widehat{f}_h - f$. In order to do this, we must recall
observation (\ref{Eq:Bias}) on the bias of $\widehat{f}_h$, and the fact
that $\|\widehat{f}_h - \mathbb{E}(\widehat{f}_h)\|_\infty=
O_{\mathrm{a.s.}}(\frac{\sqrt{\log1/h}}{(nh)^{1/2}})$ from an
application of Corollary 2.2\vspace*{-1pt} of Gin\'e and Guillou (\citeyear{GG02}). It follows
that $\|\widehat{f}_h - f\|_{\infty} \stackrel{\mathrm{a.s.}}{\rightarrow} 0$.
Similarly, $\|\mathbb{E}(\widehat{f}_h') - f'\|_{I_\delta,\infty} =
O(h^2)$, and a further application of Corollary 2.2 of Gin\'e and
Guillou (\citeyear{GG02}) gives $\|\widehat{f}_h' -
\mathbb{E}(\widehat{f}_h')\|_{I_\delta,\infty} =
O_{\mathrm{a.s.}}(\frac{\sqrt{\log1/h}}{(nh^3)^{1/2}})$. Thus
$\|\widehat{f}_h' - f'\|_{I_\delta,\infty} \stackrel{\mathrm{a.s.}}{\rightarrow}
0$. This in turn implies that with
probability one, for $n$ sufficiently large, $\widehat{f}_{h,\tau }$ is
the unique solution to $\psi(\widehat{f}_h,\widehat{f}_{h,\tau }) = 1 -
\tau $, or equivalently $\int\widehat{f}_h(x)
\mathbbm{1}_{\{\widehat{f}_h(x) \geq\widehat{f}_{h,\tau }\}} \,dx =
1-\tau $, as claimed in Section \ref{sec:asyRisk}. It remains to
note that
\[
\frac{\sum_{j=1}^r \int_{x_{2j-1}-\delta}^{x_{2j}+\delta} |\widehat
{f}_h(x) -
f(x)| \,dx}{\sum_{j=1}^{2r} |\widehat{f}_h(x_j) - f(x_j)|} = O_p(1)
\]
and
\[
\frac{\sum_{j=1}^{2r} |\widehat{f}_h(x_j) - f(x_j)|}{\min_j |\widehat
{f}_h(x_j) - f(x_j)|} = O_p(1).
\]

It follows that we can now substitute $g = \widehat{f}_h - f$ in
(\ref{Eq:PutTogether}) to deduce that
%
%
\begin{eqnarray}
\label{Eq:ProbBound}
\widehat{f}_{h,\tau } - f_{\tau } &=& \Biggl\{\sum_{j=1}^{2r}
\frac{1}{|f'(x_j)|} \Biggr\}^{-1} \Biggl\{\sum_{j=1}^{2r}
\frac{\widehat{f}_h(x_j) - f(x_j)}{|f'(x_j)|}\nonumber\\
&&\hspace*{84.4pt}{} +
\frac{1}{f_{\tau }}\sum_{j=1}^r \int_{x_{2j-1}}^{x_{2j}}
\widehat{f}_h(x) - f(x) \,dx \Biggr\} \\
&&{} + O_p \biggl(\frac{\sqrt{\log(1/h)}}{nh^2} + \frac{h^{1/2}\sqrt{\log
(1/h)}}{n^{1/2}} + h^4 \biggr).\nonumber
\end{eqnarray}
Equation (\ref{Eq:ProbBound}) shows that we can write the difference
$\widehat{f}_{h,\tau } - f_\tau $ as a sample mean of
independent and
identically distributed random variables and a small additional
remainder term. Notice from the bandwidth condition on $h^-$ in
(A2) that $\frac{\sqrt{\log(1/h)}}{nh^2} = o(h^2)$. Next,
observe that
\begin{eqnarray*}
&&\sum_{j=1}^{2r} \frac{\mathbb{E}\{\widehat{f}_h(x_j)\} -
f(x_j)}{|f'(x_j)|}
+ \frac{1}{f_{\tau }}\sum_{j=1}^r \int_{x_{2j-1}}^{x_{2j}}
\mathbb{E}\{\widehat{f}_h(x)\} - f(x) \,dx\\
&&\qquad = D_1\sum_{j=1}^{2r}
\frac{1}{|f'(x_j)|}h^2 + o(h^2),
\end{eqnarray*}
where $D_1$ is given in (\ref{Eq:D1D2D3}). Thus, in order to prove
that
%
%
\begin{equation}
\label{Eq:ExpExpansion}
\mathbb{E}(\widehat{f}_{h,\tau }) = f_{\tau } + D_1h^2 + o(h^2),
\end{equation}
it suffices by (\ref{Eq:PutTogether}) and Step \ref{Step1} to show that
for any $\eta> 0$,
\[
\mathbb{E}\bigl(|\widehat{f}_{h,\tau } - f_\tau -
D_1h^2|\mathbbm{1}_{\{\sum_{j=1}^{2r} |\widehat{f}_h(x_j) - f(x_j)| +
\|\widehat{f}_h' - f'\|_{I_\delta,\infty} > \eta\}}\bigr) = o(h^2).
\]
But this follows by Cauchy--Schwarz, because Step \ref{Step1} may be
used to show that $\mathbb{E}(\widehat{f}_{h,\tau }^2) = O(1)$, and also
\[
\mathbb{P} \Biggl(\sum_{j=1}^{2r} |\widehat{f}_h(x_j) - f(x_j)| >
\eta/2 \Biggr) + \mathbb{P}(\|\widehat{f}_h' - f'\|_{I_\delta,\infty} > \eta
/2) =
o(n^{-1}).
\]
We therefore deduce (\ref{Eq:ExpExpansion}).

In a very similar way, we can also use (\ref{Eq:PutTogether}) and the
fact that
\[
\sum_{j=1}^{2r} \frac{\operatorname{Var}\{\widehat{f}_h(x_j)\}}{f'(x_j)^2} =
\frac{D_2}{nh} \Biggl\{ \sum_{j=1}^{2r} \frac{1}{|f'(x_j)|} \Biggr\}^2 +
o \biggl(\frac{1}{nh} \biggr),
\]
where $D_2$ is given in (\ref{Eq:D1D2D3}), to deduce that
%
%
\begin{equation}
\label{Eq:Varexpansion}
\operatorname{Var}(\widehat{f}_{h,\tau }) = \frac{D_2}{nh} + o \biggl(\frac
{1}{nh} \biggr).
\end{equation}
\end{Step}
\begin{Step}\label{Step5}
We can use the results of Step \ref{Step4} to shrink
the region of interest still further. From the result of
Step \ref{Step3} we can write
%
%
\begin{eqnarray}
\label{Eq:Continuation}\qquad
&&\mathbb{E}\{\mu_f(\widehat{R}_{h,\tau } \triangle R_\tau )\} \nonumber\\
&&\qquad=
\sum_{j=1}^r \int_{x_{2j-1} - \delta_n}^{x_{2j-1} + \delta_n} f(x)
\bigl|\mathbb{P}\bigl(\widehat{f}_h(x) < \widehat{f}_{h,\tau }\bigr) -
\mathbbm{1}_{\{x < x_{2j-1}\}}\bigr| \,dx \nonumber\\
&&\qquad\quad{} + \sum_{j=1}^r \int_{x_{2j}-\delta_n}^{x_{2j}+\delta_n}
f(x) \bigl|\mathbb{P}\bigl(\widehat{f}_h(x) < \widehat{f}_{h,\tau }\bigr) -
\mathbbm{1}_{\{x \geq x_{2j}\}}\bigr| \,dx
+ o(n^{-1}) \\
&&\qquad= \frac{f_{\tau }}{(nh)^{1/2}} \sum_{j=1}^r \int_{-(nh)^{1/2}
\delta_n}^{(nh)^{1/2}\delta_n} \bigl|\mathbb{P}\bigl(\widehat{f}_h\bigl(x_{2j-1}
+ (nh)^{-1/2}t\bigr) < \widehat{f}_{h,\tau }\bigr) - \mathbbm{1}_{\{t < 0\}}\bigr|
\nonumber\\
&&\qquad\quad{} + \bigl|\mathbb{P}\bigl(\widehat{f}_h\bigl(x_{2j} +
(nh)^{-1/2}t\bigr) < \widehat{f}_{h,\tau }\bigr) - \mathbbm{1}_{\{t \geq
0\}}\bigr|\, dt + o(n^{-1}).\nonumber
\end{eqnarray}
For brevity, we write $x_j^t = x_j + (nh)^{-1/2}t$. Now, for each
$j=1,\ldots,2r$, we see that for $n$ sufficiently large,
$\mathbb{E}\{\widehat{f}_h(x_j^t) - \widehat{f}_{h,\tau }\}$ is a strictly
monotone function of $t \in[-(nh)^{1/2}\delta_n,(nh)^{1/2}\delta_n]$,
with a unique zero $t_j^*$, say. Moreover,
\[
t_j^* = \bigl\{D_1 -
\tfrac{1}{2}\mu_2(K)f''(x_j) \bigr\}\{f'(x_j)\}^{-1}
n^{1/2}h^{5/2}\{1 + o(1)\}.
\]
Fix a sequence $(t_n)$ diverging to infinity and let $I_j^n =
[-(nh)^{1/2}\delta_n,(nh)^{1/2}\delta_n] \setminus[t_j^* - t_n, t_j^*
+ t_n]$. We claim that
%
%
\begin{eqnarray}
\label{Eq:Ijn}
&&\sum_{j=1}^r \biggl\{\int_{I_{2j-1}^n}
\bigl|\mathbb{P}\bigl(\widehat{f}_h(x_{2j-1}^t)
< \widehat{f}_{h,\tau }\bigr) - \mathbbm{1}_{\{t < 0\}}\bigr|
\,dt\nonumber\\[-8pt]\\[-8pt]
&&\hspace*{26.6pt}{} +
\int_{I_{2j}^n} \bigl|\mathbb{P}\bigl(\widehat{f}_h(x_{2j}^t)
< \widehat{f}_{h,\tau }\bigr) - \mathbbm{1}_{\{t \geq0\}}\bigr|
\,dt \biggr\} \rightarrow0\nonumber
\end{eqnarray}
as $n \rightarrow\infty$. Now there exists $c_3 > 0$ such that for
all $t \in\bigcup_{j=1}^{2r} I_j^n$ and $n$ sufficiently large, we have
$|\mathbb{E}\{\widehat{f}_h(x_j^t) - \widehat{f}_{h,\tau }\}| \geq
c_3(nh)^{-1/2}t_n$. Thus there exists $c_4 > 0$ such that for all $n$
sufficiently large,
\begin{eqnarray*}
&&\bigl|\mathbb{P}\bigl(\widehat{f}_h(x_{2j-1}^t) < \widehat{f}_{h,\tau }\bigr)
- \mathbbm{1}_{\{t < 0\}}\bigr| \\[-0.8pt]
&&\qquad\leq\mathbb{P} \biggl( \biggl|\frac{\widehat{f}_h(x_{2j-1}^t)
- \mathbb{E}\{\widehat{f}_h(x_{2j-1}^t)\}}{\operatorname{Var}^{1/2}
\{\widehat{f}_h(x_{2j-1}^t)\}} \biggr| \geq c_4 t_n \biggr)\\[-0.8pt]
&&\qquad\quad{} +
\mathbb{P} \biggl( \biggl|\frac{\widehat{f}_{h,\tau } -
\mathbb{E}(\widehat{f}_{h,\tau })}{\operatorname{Var}^{1/2}
(\widehat{f}_{h,\tau })} \biggr| \geq c_4t_n \biggr) \rightarrow0,
\end{eqnarray*}
uniformly for $t \in\bigcup_{j=1}^r I_{2j-1}^n$. Since also
$|\mathbb{P}(\widehat{f}_h(x_{2j}^t) < \widehat{f}_{h,\tau })
- \mathbbm{1}_{\{t \geq0\}}| \rightarrow0$ uniformly for $t \in
\bigcup_{j=1}^{r} I_{2j}^n$, we deduce (\ref{Eq:Ijn}).
\end{Step}
\begin{Step}\label{Step6}
We also require an asymptotic expansion for
$\operatorname{Cov}(\widehat{f}_h(x_j^t),\widehat{f}_{h,\tau })$, for $t \in[t_j^*
- t_n,t_j^* + t_n]$. In fact, provided $(t_n)$ diverges sufficiently
slowly, we have
\[
\operatorname{Cov}(\widehat{f}_h(x_j^t),\widehat{f}_{h,\tau }) =
\frac{D_{3,j}}{nh} + o \biggl(\frac{1}{nh} \biggr),
\]
uniformly for $t \in[t_j^* - t_n,t_j^* + t_n]$, where $D_{3,j}$
is given at (\ref{Eq:D1D2D3}). This follows from the
expansion (\ref{Eq:PutTogether})
and the fact that provided $(t_n)$ diverges sufficiently slowly,
\begin{eqnarray*}
&&\mathbb{E} \biggl\{\frac{1}{h^2}K \biggl(\frac{x_j - X_1}{h} \biggr)
K \biggl(\frac{x_j^t - X_1}{h} \biggr) \biggr\} \\[-0.8pt]
&&\qquad= \frac{1}{h}
\int_{-\infty}^\infty K(z)K \biggl(\frac{(nh)^{-1/2}t + hz}{h} \biggr)f(x_j - hz)
\,dz \\[-0.8pt]
&&\qquad= \frac{1}{h}f_{\tau }R(K) + o(h^{-1}),
\end{eqnarray*}
uniformly for $t \in[t_j^* - t_n,t_j^* + t_n]$.
\end{Step}
\begin{Step}\label{Step7}
To complete the proof of Theorem \ref{Thm:Mainthm},
it suffices by (\ref{Eq:Continuation}) and (\ref{Eq:Ijn}) to show that
there exists a sequence $(t_n)$ diverging to infinity such that
\begin{eqnarray*}
&&\frac{f_{\tau }}{(nh)^{1/2}} \sum_{j=1}^r \biggl\{\int_{t_{2j-1}^*
- t_n}^{t_{2j-1}^* + t_n}
\bigl|\mathbb{P}\bigl(\widehat{f}_h(x_{2j-1}^t) < \widehat{f}_{h,\tau }\bigr) -
\mathbbm{1}_{\{t < 0\}}\bigr|\, dt\\[-0.8pt]
&&\hspace*{64.4pt}{} + \int_{t_{2j}^* - t_n}^{t_{2j}^* +
t_n} \bigl|\mathbb{P}\bigl(\widehat{f}_h(x_{2j}^t) <
\widehat{f}_{h,\tau }\bigr)
- \mathbbm{1}_{\{t \geq0\}}\bigr| \,dt \biggr\} \\[-0.8pt]
&&\qquad= \sum_{j=1}^{2r} \biggl[\frac{B_{1,j}\phi(B_{2,j}
n^{1/2}h^{5/2})}{(nh)^{1/2}} + B_{3,j} h^2 \{2\Phi(B_{2,j}
n^{1/2}h^{5/2}) - 1\} \biggr]\\[-0.8pt]
&&\qquad\quad{} + o \biggl(\frac{1}{(nh)^{1/2}} +
h^2 \biggr).
\end{eqnarray*}
For $i=1,\ldots,n$, let $Z_{ni}(x) =
h^{-1}K(\frac{x-X_i}{h})$ and let $\bar{Y}_n =
n^{-1}\sum_{i=1}^n Y_{ni}$, where
\begin{eqnarray*}
Y_{ni} &=& Z_{ni}(x_j^t) - f_\tau -
\Biggl\{\sum_{k=1}^{2r}\frac{1}{|f'(x_k)|} \Biggr\}^{-1}
\Biggl[\sum_{k=1}^{2r} \frac{Z_{ni}(x_k) - f(x_k)}{|f'(x_k)|}\\
&&\hspace*{154.6pt}{} +
\frac{1}{f_{\tau }}\sum_{k=1}^r \int_{x_{2k-1}}^{x_{2k}}
Z_{ni}(x) - f(x) \,dx \Biggr].
\end{eqnarray*}
By (\ref{Eq:PutTogether}) and (\ref{Eq:ProbBound}), we can write
$\widehat{f}_h(x_j^t) - \widehat{f}_{h,\tau } = \bar{Y}_n + R_n$, where $R_n
- \mathbb{E}(R_n) = o_p\{(nh)^{-1/2}\}$. Since
$\operatorname{Var}(\bar{Y}_n) = O\{(nh)^{-1}\}$ uniformly for $t \in[t_j^*
- t_n,t_j^* + t_n]$, we choose $(t_n)$ to diverge to infinity so
slowly that:
\begin{itemize}
\item$\mathbb{P} (\frac{|R_n -
\mathbb{E}(R_n)|}{\operatorname{Var}^{1/2}(\bar{Y}_n)} >
\frac{1}{t_n^2} ) \leq\frac{1}{t_n^2}$, uniformly for $t \in[t_j^*
- t_n,t_j^* + t_n]$;
\item$(nh)\operatorname{Var}(\bar{Y}_n) = R(K)f_{\tau } -2D_{3,j} +
D_2 +
o(t_n^{-1})$, uniformly for $t \in[t_j^*
- t_n,t_j^* + t_n]$;
\item$\mathbb{E}(\bar{Y}_n + R_n) = \{(nh)^{-1/2}tf'(x_j) +
D_4h^2\}\{1 + o(t_n^{-1})\}$, uniformly for $t \in[t_j^* -
t_n,t_j^* + t_n]$, where $D_4 = \frac{1}{2}\mu_2(K)f''(x_j)-D_1$;
\item$t_n =o(n^{1/6})$.
\end{itemize}
Then
\begin{eqnarray*}
&&\mathbb{P}\bigl(\widehat{f}_h(x_j^t) < \widehat{f}_{h,\tau }\bigr) -
\Phi\biggl(\frac{-tf'(x_j) - D_4n^{1/2}h^{5/2}}{\{R(K)f_{\tau } -
2D_{3,j} + D_2\}^{1/2}} \biggr) \\
&&\qquad\leq\mathbb{P} \biggl(\frac{|R_n -
\mathbb{E}(R_n)|}{\operatorname{Var}^{1/2} (\bar{Y}_n)}
> \frac{1}{t_n^2} \biggr)\\
&&\qquad\quad{} + \mathbb{P} \biggl(\frac{\bar{Y}_n
- \mathbb{E}(\bar{Y}_n)}{\operatorname{Var}^{1/2} (\bar{Y}_n)} \leq
\frac{-\mathbb{E}(\bar{Y}_n + R_n)}{\operatorname{Var}^{1/2}
(\bar{Y}_n)} + \frac{1}{t_n^2} \biggr) \\
&&\qquad\quad{} - \Phi\biggl(\frac{-tf'(x_j)
- D_4n^{1/2}h^{5/2}}{\{R(K)f_{\tau } - 2D_{3,j} + D_2\}^{1/2}} \biggr) \\
&&\qquad= O \biggl(\frac{1}{t_n^2} + \frac{1}{(nh)^{1/2}} \biggr)
+ \Phi\biggl(\frac{-\mathbb{E}(\bar{Y}_n + R_n)}{\operatorname{Var}^{1/2}
(\bar
{Y}_n)} \biggr)\\
&&\qquad\quad{}- \Phi\biggl(\frac{-tf'(x_j) - D_4n^{1/2}h^{5/2}}{\{R(K)f_{\tau }
- 2D_{3,j} + D_2\}^{1/2}} \biggr) \\
&&\qquad= o(t_n^{-1}),
\end{eqnarray*}
uniformly for $t \in[t_j^* - t_n,t_j^* + t_n]$. Here we have used
the Berry--Esseen inequality to reach the penultimate line. A very
similar argument yields a lower bound of the same order. The proof of
Step \ref{Step7}, and hence the proof of Theorem \ref{Thm:Mainthm}, is
now completed by the observation that
\begin{eqnarray*}
&&\frac{f_{\tau }}{(nh)^{1/2}} \sum_{j=1}^r
\biggl\{\int_{-\infty}^\infty
\biggl|\Phi\biggl(\frac{-tf'(x_{2j-1}) -D_4n^{1/2}h^{5/2}}{
\{R(K)f_{\tau } -2D_{3,j} + D_2\}^{1/2}} \biggr)
- \mathbbm{1}_{\{t < 0\}} \biggr| \\
&&\hspace*{45.7pt}\quad{} + \biggl|\Phi\biggl(\frac{-tf'(x_{2j})
-D_4n^{1/2}h^{5/2}}{\{R(K)f_{\tau } -2D_{3,j} + D_2\}^{1/2}} \biggr)
- \mathbbm{1}_{\{t \geq0\}} \biggr| \,dt \biggr\} \\
&&\qquad = \sum_{j=1}^{2r} \biggl[\frac{B_{1,j}
\phi(B_{2,j} n^{1/2}h^{5/2})}{(nh)^{1/2}} + B_{3,j} h^2
\{2\Phi(B_{2,j} n^{1/2}h^{5/2}) - 1\} \biggr].
\end{eqnarray*}
\end{Step}

\subsection*{Proof of Corollary \protect\ref{Cor:Opth}}

We may restrict attention to the case where $nh^5$ is bounded away
from zero and infinity. The important point to note is that under the
hypotheses of the corollary, $B_{1,j}$, $B_{2,j}$ and $B_{3,j}$ do not
depend on $j$, so we write them as $B_1$, $B_2$ and $B_3$,
respectively.

By making the substitution $x = B_2 n^{1/2}h^{5/2}$, there exist
positive constants $a = 2B_1B_2^{1/5}$ and $b=B_3/(B_1B_2)$ such that
$\lim_{n \rightarrow\infty}
n^{2/5}\mathbb{E}\{\mu_f(\widehat{R}_{h,\tau } \triangle
R_\tau )\}
= a
u(x)$, where $u(x) = x^{-1/5}\phi(x) + b
x^{4/5}\{2\Phi(x)-1\}$. Since $u$ is continuous with $u(x)
\rightarrow\infty$ as $x \searrow0$ and $x \rightarrow\infty$, it
attains its minimum in $(0,\infty)$. To show this minimum is unique,
it suffices to show that $v(x)$ has a unique zero in $(0,\infty)$,
where
\[
v(x) = \frac{5x^{6/5}}{\phi(x)} u'(x) = -1 + \frac{4bx\{2\Phi(x)-1\}
}{\phi(x)} + 5(2b-1)x^2.
\]
Now we have
\begin{eqnarray*}
v'(x) &=& 2(14b-5)x + \frac{4b(1+x^2)\{2\Phi(x)-1\}}{\phi(x)}, \\
v''(x) &=& 2(18b-5) + 8bx^2 + \frac{4b(3x+x^3)\{2\Phi(x)-1\}}{\phi(x)}.
\end{eqnarray*}
There are therefore two cases to consider: if $b \geq5/18$, then $v$
is strictly convex, so since $v(0+)=-1$ and $v(x) \rightarrow\infty$
as $x \rightarrow\infty$, we see that $v$ has a unique zero in
$(0,\infty)$. On the other hand, if $b < 5/18$, then there exists
$x^* \in(0,\infty)$ such that $v''(x) < 0$ for $x \in(0,x^*)$ and
$v''(x) > 0$ for $x \in(x^*,\infty)$. But if $b < 5/18$ then $v'(x)
< 0$, for sufficiently small $x > 0$, so from $v(0+) = -1$, it again
follows that $v$ has a unique zero.

Write $x_{\min}$ for the unique minimum of $u$ in
$(0,\infty)$, and let $c_{\mathrm{opt}} =
(x_{\min}/B_2)^{2/5}$. We conclude that any optimal bandwidth
sequence $(h_{\mathrm{opt}})$, in the sense of minimizing
$\mathbb{E}\{\mu_f(\widehat{R}_{h,\tau } \triangle R_\tau )\}
$, must
satisfy $h_{\mathrm{opt}} = c_{\mathrm{opt}}n^{-1/5}\{1 + o(1)\}$
as $n \rightarrow\infty$.

\subsection*{Proof of Theorem \protect\ref{Thm:Hath}}

We require a bound on $|\widehat{x}_{j,h_0} - x_j|$ for $j=1,\ldots,2r$.
To this end, let $\tilde{f} = f + g$ be another density satisfying the
same conditions as~$f$. From Step \ref{Step4} of the proof of
Theorem \ref{Thm:Mainthm}, we see that for sufficiently small values
of $\|g\|_\infty+ \|g'\|_{I_\delta,\infty}$, there exist precisely
$2r$ values $\tilde{x}_1 < \cdots< \tilde{x}_{2r}$ such that
$\tilde{f}(\tilde{x}_j) = \tilde{f}_\tau $. Moreover, provided
$\sum_{j=1}^r \int_{x_{2j-1}-\delta}^{x_{2j}+\delta} |g(x)| \,dx =
O(\sum_{j=1}^{2r} |g(x_j)|)$ as $\sum_{j=1}^{2r} |g(x_j)| +
\|g'\|_{I_\delta,\infty} \rightarrow0$, we have $\tilde{x}_j - x_j =
O(|g(x_j)|)$ as $\sum_{j=1}^{2r} |g(x_j)| + \|g'\|_{I_\delta,\infty}
\rightarrow0$. Substituting $\tilde{f} = \widehat{f}_{h_0}$, so that $g =
\widehat{f}_{h_0} - f$ and $\tilde{x}_j = \widehat{x}_{j,h_0}$, we have
$|\widehat{x}_{j,h_0} - x_j| = O_p(n^{-2/5})$.

It follows that $\widehat{D}_1 = D_1 + O_p(n^{-2/9})$, the crucial fact
being that $\widehat{f}_{h_2}''(\widehat{x}_{j,h_0}) - f''(x_j) =
O_p(n^{-2/9})$. Similarly, $\widehat{D}_2 = D_2 + O_p(n^{-2/7})$ and
$\widehat{D}_{3,j} = D_{3,j} + O_p(n^{-2/7})$ for $j=1,\ldots,2r$. Thus
$\widehat{B}_{1,j} = B_{1,j} + O_p(n^{-2/7})$, $\widehat{B}_{2,j} = B_{2,j} +
O_p(n^{-2/9})$ and $\widehat{B}_{3,j} = B_{3,j} + O_p(n^{-2/9})$. We
deduce that for any $0 < c_1 <c_2 < \infty$, we have
$\widehat{\mathrm{AR}}_n(c) = \mathrm{AR}(c)\{1 + O_p(n^{-2/9})\}$, uniformly
for $c \in[c_1,c_2]$, and a standard Taylor expansion argument then
gives that $\widehat{c}_{\mathrm{opt}} = c_{\mathrm{opt}}\{1 +
O_p(n^{-2/9})\}$. Both conclusions of the theorem follow
immediately.

\subsection*{Proof of Theorem \protect\ref{Thm:HathHO}}

Let $z_n = \delta/h_2$, where $\delta$ is small enough that $f$ has 12
continuous derivatives in $\bigcup_{j=1}^{2r} [x_j-\delta,x_j+\delta]$.
Under the conditions of the theorem, we may integrate by parts twice
and apply a Taylor expansion to obtain
\begin{eqnarray*}
|\mathbb{E}\{\widehat{f}_{h_2}''(x_j)\} -f''(x_j)|
&=& \biggl|\int_{-z_n}^{z_n} K_2(z)
\{f''(x_j-h_2z) - f''(x_j)\} \,dz \biggr| + o(h_2^{10})\\
&=& O(h_2^{10}).
\end{eqnarray*}
This expression for the bias can be combined with the standard fact
that $\operatorname{Var} \widehat{f}_{h_2}''(x_j) = O\{(nh_2^5)^{-1}\}$ and the
bound on $|\widehat{x}_{j,h_0} - x_j|$ from the proof of
Theorem \ref{Thm:Hath} to yield $\widehat{f}_{h_2}''(\widehat{x}_{j,h_0}) -
f''(x_j) = O_p(n^{-2/5})$. Similar computations give
$\widehat{f}_{h_1}'(\widehat{x}_{j,h_0}) - f'(x_j) = O_p(n^{-2/5})$. The rest
of the proof mirrors the proof of Theorem \ref{Thm:Hath}.
\end{appendix}

\section*{Acknowledgments}

The authors are grateful to Tarn Duong, Inge Koch,
Steve Marron and Richard Nickl for their comments on aspects of this
research, and to the organizers of a workshop on statistical
research held at the Keystone Resort, Colorado, USA, on 4th--8th June, 2007.

\printaddresses

\end{document}